\documentclass[11pt, cjk, twoside, leqno]{article}

\usepackage{stmaryrd}
\usepackage{amsfonts}
\usepackage{amssymb}
\usepackage{amsmath}
\usepackage{amsthm}
\usepackage{xcolor}

\usepackage{mathrsfs}

\usepackage{CJK}

\allowdisplaybreaks

\def\rr{{\mathbb R}}

\def\zz{{\mathbb Z}}
\def\cc{{\mathbb C}}

\def\nn{{\mathbb N}}

\def\cj{{\mathcal J}}

\def\scs{{\mathscr S}}
\def\scp{{\mathscr P}}

\def\fz{\infty}
\def\az{\alpha}

\def\lz{\lambda}
\def\blz{\Lambda}

\def\vz{\varphi}

\def\lf{\left}
\def\r{\right}

\def\hs{\hspace{0.25cm}}
\def\ls{\lesssim}

\def\noz{\nonumber}
\def\wz{\widetilde}

\def\st{\subset}

\def\bh{\backslash}

\def\gfz{\genfrac{}{}{0pt}{}}

\def\supp{\mathop\mathrm{\,supp\,}}

\def\pd{{\scp(\rr^D)}}
\def\sch{{\scs(\rr^D)}}
\def\schd{{\scs'(\rr^D)}}
\def\schi{{\scs_{\infty}(\rr^D)}}
\def\schid{{\scs'_{\infty}(\rr^D)}}

\def\lon{{L^1(\rr^D)}}
\def\ltw{{L^2(\rr^D)}}
\def\li{{L^{\infty}(\rr^D)}}
\def\bmo{\mathop\mathrm{\,{\rm BMO}(\rr^D)}}

\def\hon{{H^1(\rr^D)}}

\def\loq{{\mathop\mathrm{WE}^{1,\,q}(\rr^D)}}
\def\liq{{\mathop\mathrm{WE}^{\infty,\,q}(\rr^D)}}

\def\tls{{\dot{F}^0_{1,\,q}(\rr^D)}}
\def\dtl{{\dot{F}^0_{\fz,\,q'}(\rr^D)}}

\newtheorem{theorem}{Theorem}[section]
\newtheorem{lemma}[theorem]{Lemma}

\newtheorem{proposition}[theorem]{Proposition}
\theoremstyle{definition}
\newtheorem{remark}[theorem]{Remark}
\newtheorem{definition}[theorem]{Definition}
\newtheorem{example}[theorem]{Example}

\numberwithin{equation}{section}

\begin{document}

\arraycolsep=1pt

\title{\vspace{-2cm}\Large\bf
The duality about function set and Fefferman-Stein
Decomposition \footnotetext {\hspace{-0.35cm}
2010 {\it Mathematics Subject Classification}. Primary 42B35;
Secondary 46E35, 42B20.
\endgraf {\it Key words and phrases}.
Riesz transform, Fefferman-Stein
decomposition, Triebel-Lizorkin space.
\endgraf
This project is supported by the National
Natural Science Foundation of China
(Grant No. 11571261) and Macao Government FDCT099.}}
\author{
Qixiang Yang \ \ and \ \ Tao Qian \,\footnote{Corresponding author}
 }
\date{  }
\maketitle

\vspace{-0.8cm}

\begin{center}
\begin{minipage}{13cm}
{\small {\bf Abstract}\quad Let $D\in\mathbb{N}$, $q\in[2,\infty)$ and
$(\mathbb{R}^D,|\cdot|,dx)$ be the Euclidean space equipped
with the $D$-dimensional Lebesgue measure. In this article,
the authors establish the Fefferman-Stein
decomposition of Triebel-Lizorkin spaces $\dot{F}^0_{\infty,\,q'}(\mathbb{R}^D)$ on basis of the dual on function set which has special topological structure.
The function in Triebel-Lizorkin spaces $\dot{F}^0_{\infty,\,q'}(\mathbb{R}^D)$ can be written as the certain combination of $D+1$ functions in $\dot{F}^0_{\infty,\,q'}(\mathbb{R}^D) \bigcap L^{\infty}(\mathbb{R}^D)$.

To get such decomposition, {\bf (i),}
The authors introduce some auxiliary function space $\mathrm{WE}^{1,\,q}(\mathbb R^D)$ and
$\mathrm{WE}^{\infty,\,q'}(\mathbb{R}^D)$ defined via wavelet expansions.
The authors proved $\tls\subsetneqq L^{1}(\rr^D)  \bigcup \tls \subset {\rm WE}^{1,\,q}(\rr^D)\subset L^{1}(\rr^D) +    \tls$ and $\mathrm{WE}^{\infty,\,q'}(\mathbb{R}^D)$ is strictly contained in
$\dot{F}^0_{\infty,\,q'}(\mathbb{R}^D)$.
{\bf (ii),} The authors establish the Riesz transform characterization
of Triebel-Lizorkin spaces $\dot{F}^0_{1,\,q}(\mathbb{R}^D)$ by function set $\mathrm{WE}^{1,\,q}(\mathbb R^D)$.
{\bf (iii),} We also consider  the dual of $\mathrm{WE}^{1,\,q}(\mathbb R^D)$.
As a consequence of the above results, the authors get also Riesz transform characterization
of Triebel-Lizorkin spaces $\dot{F}^0_{1,\,q}(\mathbb{R}^D)$
by Banach space $L^{1}(\rr^D) +    \tls$.

Although Fefferman-Stein type decomposition
when $D=1$ was obtained by
C.-C. Lin et al. [Michigan Math. J. 62 (2013), 691-703], as was pointed
out by C.-C. Lin et al., the approach used in the
case $D=1$ can not be applied to the cases $D\ge2$, which needs some new methodology. 
}
\end{minipage}
\end{center}

\section{Introduction and main results}\label{s1}

\hskip\parindent The Riesz transforms on $\rr^D$ ($D\ge2$), which are natural generalizations
of the Hilbert transform on $\rr$, may be the
most typical examples of Calder\'on-Zygmund operators
(see, for example, \cite{g08,s70,s93} and references therein). It
is well known that the Riesz transforms have many interesting properties,
for example, they are the simplest, nontrivial, ¡°invariant¡± operators
under the action of the group of rotations in the Euclidean space $\rr^D$,
and they also constitute typical and important examples of Fourier multipliers.
Moreover, they can be used to mediate between various combinations
of partial derivatives of functions. All these
properties make the Riesz transforms ubiquitous in mathematics and useful in various fields of analysis
such as partial differential equations and harmonic analysis
(see \cite{s70,s93} for more details on their applications).

The Riesz transform characterization of Hardy spaces plays
important roles in the real variable theory of Hardy spaces (see,
for example, \cite{fs, s70}). Via this Riesz transform characterization of
the Hardy space $\hon$ and the duality between $\hon$ and
the space of functions with bounded mean oscillation, $\bmo$,
Fefferman and Stein \cite{fs} further obtained the nowadays so-called Fefferman-Stein
decomposition of $\bmo$. Later, Uchiyama \cite{u82} gave a constructive proof
of the Fefferman-Stein decomposition of $\bmo$. Since then,
many articles focus on the classical Riesz transform characterization
and the Fefferman-Stein decomposition of different variants of Hardy spaces
and BMO spaces; see, for example, \cite{ccyy,cg,g79,ll,yzn} and references therein.
Recently, Lin et al. \cite{lly} established the Hilbert transform
characterization of Triebel-Lizorkin spaces ${\dot{F}^0_{1,\,q}(\rr)}$
and the Fefferman-Stein decomposition of Triebel-Lizorkin spaces
${\dot{F}^0_{\fz,\,q'}(\rr)}$ for each $q\in[2,\fz)$.
Yang et al. \cite{yql} obtained the Fefferman-Stein
decomposition for $Q$-spaces $Q_\az(\rr^D)$ and
the Riesz transform characterization of $P^\az(\rr^D)$,
the predual of $Q_\az(\rr^D)$, for any $\az\in[0,\fz)$.

As was pointed out by Lin et al. in \cite[Remark 1.4]{lly}, the approach used
in \cite{lly} for the Hilbert transform characterization of Triebel-Lizorkin
spaces ${\dot{F}^0_{1,\,q}(\rr)}$
can not be applied to $\tls$ when $D\ge2$, which needs to
develop some new skills. In this article, motivated by some ideas
from \cite{lly,yql}, we establish the Riesz transform
characterization of Triebel-Lizorkin spaces ${\dot{F}^0_{1,\,q}(\rr^D)}$
and the Fefferman-Stein decomposition of Triebel-Lizorkin spaces
${\dot{F}^0_{\fz,\,q'}(\rr^D)}$ for all $D\in\nn:=\{1,\,2,\,\ldots\}$
and $q\in[2,\fz)$.

In order to state the main results of this article,
we now recall the definition
of the Triebel-Lizorkin space
$\tls$ from \cite{tr1}; see also \cite{tr2,tr3,tr4,fjw}.
Let $\sch$ and $\schd$ be the \emph{Schwartz space} and
its \emph{dual} respectively, and $\pd$ the \emph{class of all polynomials on $\rr^D$}.
Following \cite{tr1}, we also let
$$
\schi:=\lf\{\vz\in\sch:\ \int_{\rr^D}\vz(x)x^{\az}\,dx=0
\ {\rm for\ all}\ \az\in\zz^D_+\r\}
$$
and $\schid$ be its dual. Here and hereafter, $\zz_+:=\nn\cup\{0\}$,
$\zz^D_+:=(\zz_+)^D$ and, for any $\az:=(\az_1,\ldots,\az_D)\in\zz^D_+$
and $x:=(x_1,\ldots,x_D)\in\rr^D$, $x^\az:=x_1^{\az_1}\cdots x_D^{\az_D}$.

\begin{definition}\label{da.a}
Let $\vz\in\sch$ satisfy $\supp(\widehat{\vz})\st\{\xi\in\rr^D:\
\frac12 \le|\xi|\le 2\}$, $|\widehat{\vz}(\xi)|\ge c>0$ if
$\frac35 \le|\xi|\le \frac53$, and
$\sum_{j\in\zz}|\widehat{\vz}(2^j\xi)|=1$ if $\xi\neq0$,
where $c$ is a positive constant.
Write $\vz_j(\cdot):=2^{Dj}\vz(2^{j}\cdot)$ for any $j\in\zz$.
Let $q\in(1,\fz)$. Then the \emph{homogeneous Triebel-Lizorkin space}
$\tls$ is defined to be the set of all $f\in\schid$ such that
$$
\|f\|_{\tls}:=\lf\|\lf\{\sum_{j\in\zz}
\lf|\vz_j\ast f\r|^q\r\}^{1/q}\r\|_{\lon}<\fz.
$$
\end{definition}

\begin{remark}\label{ra.n}
(i) It is well known that $\schid=\schd/\pd$ with equivalent
topologies; see, for example, \cite[Proposition 8.1]{ywy}
and \cite[Theorem 6.28]{s16} for an exact proof.

(ii) From \cite[p.\,42]{fjw}, it follows that
$\dot{F}^0_{1,\,2}(\rr^D)=\hon$ with equivalent norms.
Obviously, for any $q\in[2,\fz)$, $\hon\st\tls$.
\end{remark}

Now we recall the definition of the dual space of $\tls$, $\dtl$,
from \cite[p.\,70]{fj}, where $1/q+1/q'=1$.

\begin{definition}\label{da.b}
Let $q\in(1,\fz)$. Then the \emph{homogeneous Triebel-Lizorkin space}
$\dot{F}^{0,\,q}_\fz(\rr^D)$ is defined to be the set
of all $f\in\schid$ such that
$$
\|f\|_{\dot{F}^0_{\fz,\,q}(\rr^D)}:=\sup_{\{Q:\ \rm dyadic\ cube\}}
\lf\{\frac1{|Q|}\int_{Q}\sum_{j=-\log_2\ell(Q)}^\fz\lf|\vz_j\ast f(x)\r|^q
\,dx\r\}^{1/q}<\fz,
$$
where the supremum is taken over all dyadic cubes $Q$ in $\rr^D$ and
$\ell(Q)$ denotes the \emph{side length} of $Q$.
\end{definition}

\begin{remark}\label{ra.c}
(i) From \cite[p.\,42]{fjw}, it follows that $\dot{F}^0_{\fz,\,2}(\rr^D)=\bmo$
with equivalent norms.

(ii) It was shown in \cite[(5.2)]{fj} that, for each $q\in(1,\fz)$,
$\dtl$ is the dual space of $\tls$. In particular, $\bmo$ is
the dual space of $\hon$, which was proved before in \cite{fs}.
\end{remark}

Next we recall the definition of the $1$-dimensional
Meyer wavelets from \cite{w97}; see also \cite{lly, m92, QY} for a
different version.
Let $\Phi\in C^\fz(\rr)$, the \emph{space of all infinitely differentiable
functions on $\rr$},  satisfy
\begin{equation}\label{ph1}
0\le\Phi(\xi)\le\frac1{\sqrt{2\pi}}\quad {\rm for\ any\ }\xi\in\rr,
\end{equation}
\begin{equation}\label{ph2}
\Phi(\xi)=\Phi(-\xi)\quad {\rm for\ any\ }\xi\in\rr,
\end{equation}
\begin{equation}\label{ph3}
\Phi(\xi)=\frac1{\sqrt{2\pi}}\quad {\rm for\ any\ }\xi\in[-2\pi/3,2\pi/3],
\end{equation}
\begin{equation}\label{ph4}
\Phi(\xi)=0\quad {\rm for\ any\ }\xi\in(-\fz,4\pi/3]\cup[4\pi/3,\fz),
\end{equation}
and
\begin{equation}\label{ph5}
\lf[\Phi(\xi)\r]^2+\lf[\Phi(\xi-2\pi)\r]^2=\frac1{2\pi}
\quad {\rm for\ any\ }\xi\in[0,2\pi].
\end{equation}
In what follows, the \emph{Fourier transform} and
the \emph{reverse Fourier transform}
of a suitable function $f$ on $\rr^D$ are defined by
$$
\widehat{f}(\xi):=(2\pi)^{-D/2}\int_{\rr^D}e^{-i\xi x}f(x)\,dx \quad
{\rm for\ any\ }\xi\in\rr^D,
$$
respectively,
$$
\check{f}(x):=(2\pi)^{-D/2}\int_{\rr^D}e^{ix\xi}f(\xi)\,d\xi \quad
{\rm for\ any\ }x\in\rr^D.
$$

From \cite[Proposition 3.2]{w97}, it follows that $\phi:=\check{\Phi}$
(the \emph{``farther" wavelet}) is a scaling function of a
\emph{multiresolution analysis} defined as in
\cite[Definition 2.2]{w97}.
The \emph{corresponding function} $m_{\phi}$ of $\phi$, satisfying
$\widehat{\phi}(2\cdot)=m_{\phi}(\cdot)\widehat{\phi}(\cdot)$,
is a $2\pi$-periodic function which equals $\sqrt{2\pi}\Phi(2\cdot)$
on the interval $[-\pi,\pi)$.

Furthermore, by \cite[Theorem 2.20]{w97}, we construct a
$1$-dimensional wavelet $\psi$ (the \emph{``mother" wavelet}) by
setting
$\widehat{\psi}(\xi):=e^{i\xi/2}m_{\phi}(\xi/2+\pi)\Phi(\xi/2)$
for any $\xi\in\rr$.
It was shown in \cite[Proposition 3.3]{w97} that $\psi$ is a
real-valued $C^{\fz}(\rr)$ function, $\psi(-1/2-x)=\psi(-1/2+x)$
for all $x\in\rr$, and
\begin{equation}\label{ps1}
\supp\lf(\widehat{\psi}\r)\st[-8\pi/3,-2\pi/3]\cup[2\pi/3,8\pi/3].
\end{equation}
Such a wavelet $\psi$ is called a \emph{$1$-dimensional Meyer wavelet}.

Let $D\in\nn\cap[2,\fz)$ and
$\vec{0}:=(\overbrace{0,\ldots,0}^{D\ {\rm times}})$.
The $D$-dimensional Meyer wavelets are constructed
by tensor products as follows.
Let $x:=(x_1,\ldots,x_D)\in\rr^D$,
$E_D:=\{0,1\}^D\bh\{\vec{0}\}$
and, for any $\lz:=(\lz_1,\ldots,\lz_D)\in E_D$, define
$$
\psi^{\lz}(x):=\phi^{\lz_1}(x_1)\cdots\phi^{\lz_D}(x_D),
$$
with $\phi^{\lz_j}(x_j):=\phi(x_j)$ if $\lz_j=0$ and
$\phi^{\lz_j}(x_j):=\psi(x_j)$ if $\lz_j=1$. As in \cite{w97},
for any $(\lz,j,k)\in\blz_D:=\{(\lz,j,k):\ \lz\in E_D,\ j\in\zz,\ k\in\zz^D\}$
and $x\in\rr^D$, we let $\psi^\lz_{j,\,k}(x):=2^{Dj}\psi^{\lz}(2^jx-k)$
and, for $\lz=\vec{0}$ and any $k:=(k_1,\ldots,k_D)$,
let $\psi^{\vec{0}}_{j,\,k}(x)
:=2^{Dj}\phi(2^jx_1-k_1)\cdots\phi(2^jx_D-k_D)$ and
$\psi^{\vec{0}}(x):=\phi(x_1)\cdots\phi(x_D)$.

By \cite[Proposition 3.1]{w97} and 
arguments of tensor products,
we know that, for any $(\lz,j,k)\in\blz_D$,
$\psi^{\lz}_{j,\,k}\in\schi$. Thus, for any $(\lz,j,k)\in\blz_D$
and any $f\in\schid$, let
$a^{\lz}_{j,\,k}(f):=\langle f,\psi^{\lz}_{j,\,k}\rangle$,
where $\langle \cdot,\cdot\rangle$ represents the duality
between $\schid$ and $\schi$.
From the proof of \cite[Theorem (7.20)]{fjw}, it follows that, for any
$f\in\schid$,
\begin{equation}\label{x.x}
f=\sum_{\lz\in E_D}\sum_{j\in\zz}\sum_{k\in\zz^D}
a^{\lz}_{j,\,k}(f)\psi^{\lz}_{j,\,k}\quad {\rm in}\quad
\schid.
\end{equation}
Moreover, by \cite[Proposition 5.2]{w97}, we know that
$\{\psi^{\lz}_{j,\,k}\}_{(\lz,j,k)\in\blz_D}$
is an orthonormal basis of $\ltw$.

For any $\ell\in\{1,\ldots,D\}$ and any $f\in\sch$,
denote by $R_{\ell}(f)$ the \emph{Riesz transform}
of $f$, which is defined by setting
$$
\widehat{R_\ell(f)}(\xi):=-i\frac{\xi_\ell}{|\xi|}\widehat{f}(\xi)
\quad {\rm for\ any\ }\xi\in\rr^D.
$$
Since \eqref{ph4} and \eqref{ps1} hold true, by \cite[(5.2)]{yql},
we know that, for any
$\ell\in\{1,\ldots,D\}$, $(\lz,j,k),\,(\wz{\lz},\wz{j},\wz{k})\in\blz_D$
and $|j-\wz{j}|\ge2$, we have
\begin{equation}\label{b.d}
\lf( R_{\ell}\lf(\psi^{\lz}_{j,\,k}\r),\psi^{\wz{\lz}}_{\wz{j},\,\wz{k}}
\r)=0,
\end{equation}
where $(\cdot,\cdot)$ denotes the inner product in $\ltw$.

Now we recall the wavelet characterization of $\tls$ and $\dot{F}^{0}_{\infty,q}(\mathbb{R}^{D})$
(see, for example, \cite[Theorem (7.20)]{fjw}). For $j\in \mathbb{Z}$ and $k=(k_1,\cdots, k_{D})\in \mathbb{Z}^{D}$, denote
$Q_{j,k}=\prod\limits^{D}_{l=1} [2^{-j}k_l, 2^{-j}(1+k_{l})[$.

\begin{theorem}\label{ta.d}
Let $q\in(1,\fz)$. Then

(i) $f\in\tls$ if and only if $f\in\schid$ and
$$
\cj_f:=\lf\|\lf\{\sum_{(\lz,\,j,\,k)\in\blz_D}
\lf[2^{Dj}\lf|a^{\lz}_{j,\,k}(f)\r|\chi\lf(2^jx-k\r)\r]^q
\r\}^{1/q}\r\|_{\lon}<\fz,
$$
where $\chi$ denotes the characteristic function of the cube $[0,1)^D$.
Moreover, there exists a positive constant $C$
such that, for all $f\in\tls$,
$$
\frac1C\|f\|_{\tls}\le\cj_f\le C\|f\|_{\tls}.
$$

(ii) $f\in \dot{F}^{0}_{\infty,q}(\mathbb{R}^{D})$ if and only if $f\in\schid$ and there exists $C>0$ such that for all dyadic cube $Q$,
$$ \sum_{(\lz,\,j,\,k)\in\blz_D, Q_{j,k} \subset Q} 2^{(q-1) j D} |a^{\lambda}_{j,k}|^{q} \leq C|Q|.
$$
\end{theorem}

\begin{remark}\label{ra.o}
By Remark \ref{ra.n}(ii) and Theorem \ref{ta.d}, we also obtain
the wavelet characterization of $\hon$ as in \cite[p.\,143]{m92}.
\end{remark}

To consider Fefferman-Stein type decomposition for $\dot{F}^{0}_{\infty,q}(\mathbb{R}^{D})$, we need to study some properties relative to frequency.
Hence we use Meyer wavelets to introduce
the auxiliary function spaces $\loq$. 
We consider the linear functional on these function sets and consider some exchangeability of Riesz transform and some sums of orthogonal projector operator defined by Meyer wavelets.

Let $q\in(1,\fz)$ and $f\in\schid$. For any $s\in\zz$,
$N\in\nn$ and $t\in\{0,\ldots,N+1\}$, let
\begin{equation}\label{b.a}
P_{s,\,N}f:=\sum_{\{(\lz,\,j,\,k)\in\blz_D:\ s-N\le j\le s\}}
a^{\lz}_{j,\,k}(f)\psi^{\lz}_{j,\,k}
\quad {\rm in}\quad \schid.
\end{equation}
For each $t\in\{0,\ldots,N+1\}$, let, in $\schid$,
\begin{equation}\label{b.x}
T^{(1)}_{s,\,t,\,N}(f):=\begin{cases}
0, \ \  \ \ &t=0, \\
\displaystyle\sum_{\{(\lz,\,j,\,k)\in\blz_D:\ s-t+1\le j\le s\}}
a^{\lz}_{j,\,k}(f)\psi^{\lz}_{j,\,k}, \ \  \ \
&t\in\{1,\ldots,N+1\}
\end{cases}
\end{equation}
and
\begin{equation}\label{b.y}
T^{(2)}_{s,\,t,\,N}(f):=\begin{cases}
\displaystyle\sum_{\{(\lz,\,j,\,k)\in\blz_D:\ s-N\le j\le s-t\}}
a^{\lz}_{j,\,k}(f)\psi^{\lz}_{j,\,k}, \ \  \ \ &t\in\{0,\ldots,N\}, \\
0, \ \  \ \ &t=N+1.
\end{cases}
\end{equation}

\begin{definition}\label{da.e}
Then the \emph{space $\liq$} is defined to be the space of all $f\in\schid$
such that
$$
\|f\|_{\liq}:=\sup_{\{s\in\nn,\,N\in\nn\}}
\sup_{t\in\{0,\ldots,N+1\}}
\lf[\lf\|T^{(1)}_{s,\,t,\,N}(f)\r\|_{{\dot{F}^0_{\fz,\,q}(\rr^D)}}
+\lf\|T^{(2)}_{s,\,t,\,N}(f)\r\|_{\li}\r]<\fz.
$$
\end{definition}

It is easy to see that
\begin{proposition}
For $1<q\leq \infty$, $\liq= L^{\infty}(\mathbb{R}^{D}) \bigcap \dot{F}^{0}_{\infty,q}(\mathbb{R}^{D})$ are Banach spaces.
\end{proposition}

\begin{definition}\label{da.e111}
The relative \emph{space $\loq$} is defined to be the space of all $f\in\schid$
such that
$$
\|f\|_{\loq}:=\sup_{\{s\in\nn,\,N\in\nn\}}
\min_{t\in\{0,\ldots,N+1\}}\lf[\lf\|T^{(1)}_{s,\,t,\,N}(f)\r\|_{\tls}
+\lf\|T^{(2)}_{s,\,t,\,N}(f)\r\|_{\lon}\r]<\fz.
$$
Further, for $f\in L^{1}(\mathbb{R}^{D}) \bigcup \dot{F}^{0}_{1,q}(\mathbb{R}^{D})$, we define
$$\|f\|_{\{1, q\}}:= \min (\|f\|_{L^{1}}, \|f\|_{\dot{F}^{0}_{1,q}}).$$
For $f\in L^{1}(\mathbb{R}^{D}) + \dot{F}^{0}_{ 1,q}(\mathbb{R}^{D})$, we define
$$\|f\|_{1,q}:= \inf\limits_{f+g\in L^{1} + \dot{F}^{0}_{1,q}} \{\|f\|_{L^{1}} + \|g\|_{\dot{F}^{0}_{1,q}}\}.$$
\end{definition}

\begin{remark}\label{ra.k}
The spaces $\loq$ and $\liq$, with $q\in(1,\fz)$, when
$D=1$ were introduced by Lin et al. \cite[p.\,693]{lly},
respectively, \cite[p.\,694]{lly}, which were denoted by $L^{1,\,q}(\rr)$,
respectively, $L^{\fz,\,q}(\rr)$. To distinguish these spaces with
the well-known Lorentz spaces, we use the notation
$\loq$ and $\liq$ which indicate that these spaces are defined
via wavelet expansions. Recall also that the space $\loq$
was also called the relative $L^1$ space in \cite[p.\,693]{lly}.


\end{remark}

We know,
$\forall N\geq 1$, the function $P_{N}f(x)= \sum\limits_{(\epsilon,j,k)\in \Lambda_{D}, |j|+|k|\leq 2^{N}} a^{\lambda}_{j,k} \psi^{\lambda}_{j,k}(x)\in \schi$.
Set $ A= \loq$ or $ L^{1}(\mathbb{R}^{D}) \bigcup \dot{F}^{0}_{1,q}(\mathbb{R}^{D})$ or $L^{1}(\mathbb{R}^{D}) + \dot{F}^{0}_{1,q}(\mathbb{R}^{D})$.
If $f\in A$, then $P_{N}f\in A$.
It is easy to see that
\begin{proposition} \label{pr.1.9} For $1\leq q<\infty$,

(i) $\loq$ is complete with the above induced norm.

(ii) The functions in $\schi$ are dense in $A$.
\end{proposition}

\begin{remark}\label{re:Banach,function}
Let $q\in[2,\fz)$. It was shown in \cite[p.\,239]{tr3} that the dual space of
$\tls$ is   $\dtl$. Further $\dot{F}^{0}_{1,2}(\mathbb{R}^{D})= H^{1}(\mathbb{R}^{D}) \subset L^{1}(\mathbb{R}^{D}).$
Hence \begin{equation}\label{eq:Hardy}
L^{1}(\mathbb{R}^{D})\bigcup \dot{F}^{0}_{1,2}(\mathbb{R}^{D})={\rm WE}^{1,2} (\mathbb{R}^{D})= L^{1}(\mathbb{R}^{D})+ \dot{F}^{0}_{1,2}(\mathbb{R}^{D})=L^{1}(\mathbb{R}^{D}).\end{equation}
Let $2<q<\infty$. 
$L^{1}(\mathbb{R}^{D})+ \dot{F}^{0}_{1,q}(\mathbb{R}^{D})$ are Banach spaces.
$WE^{1,q} (\mathbb{R}^{D})$ are function sets, not Banach spaces. Moreover,
the following equalities are {\bf not} true
$$L^{1}(\mathbb{R}^{D})\bigcup \dot{F}^{0}_{1,q}(\mathbb{R}^{D})= {\rm WE}^{1,q} (\mathbb{R}^{D})= L^{1}(\mathbb{R}^{D})+ \dot{F}^{0}_{1,q}(\mathbb{R}^{D}).$$
In fact, the above two equal signs both have to be changed to the inclusion sign ``$\subset$".
\end{remark}


For $A$, we can use distributions to define their dual elements.
\begin{definition}
For $1\leq q<\infty$ and the function set $A\subset L^{1}(\mathbb{R}^{D}) + \dot{F}^{0}_{1,q}(\mathbb{R}^{D})$,  we call
$l$ to be a dual element of $A$,
if $l\in \schid$ and
$$\sup\limits_{f\in \schi, \|f\|_{A}\leq 1} |\langle l, f\rangle|<\infty.$$
We write $l\in A'$.
\end{definition}
$A'$ is a linear space. In fact, for $\alpha, \beta\in \mathcal{C}$ and $l_1, l_2\in A'$, we know that $\alpha l_1 +\beta l_2\in A'$.
Further, $L^{1}(\mathbb{R}^{D})+ \dot{F}^{0}_{1,q}(\mathbb{R}^{D})$ is the linearization function space of the set
$L^{1}(\mathbb{R}^{D})\bigcup \dot{F}^{0}_{1,q}(\mathbb{R}^{D})$ or the set $ {\rm WE}^{1,q} (\mathbb{R}^{D}).$
The dual elements on the set $L^{1}(\mathbb{R}^{D})\bigcup \dot{F}^{0}_{1,q}(\mathbb{R}^{D})$ or on the set
${\rm WE}^{1,q} (\mathbb{R}^{D})$ are the same as which on the linear space $L^{1}(\mathbb{R}^{D})+ \dot{F}^{0}_{1,q}(\mathbb{R}^{D})$.
Now we are ready to state the first main auxiliary result of this paper.

\begin{theorem}\label{th:111}
For $q\in [2,\infty)$, we have
$$\big(L^{1}(\mathbb{R}^{D})\bigcup \dot{F}^{0}_{1,q}(\mathbb{R}^{D})\big)'= \big({\rm WE}^{1,q} (\mathbb{R}^{D})\big)'=\big( L^{1}(\mathbb{R}^{D})+ \dot{F}^{0}_{1,q}(\mathbb{R}^{D})\big)'= L^{\infty}(\mathbb{R}^{D}) \bigcap \dot{F}^{0}_{\infty,q'}(\mathbb{R}^{D}).$$

\end{theorem}

For $q=2$, due to the equation (\ref{eq:Hardy}), the above Theorem \ref{th:111} is evident. For general $q$, the proof of this theorem will be given in the final section.

Next we state the second main auxiliary result which will be needed in the proof of our Fefferman-Stein type decomposition.
We will use certain exchangeability of Meyer wavelets and Riesz transform to prove Theorem \ref{ta.h} in section 2.

\begin{theorem}\label{ta.h}
Let $D\in\nn$ and $q\in[2,\fz)$.
Then $f\in\schid$ belongs to $\tls$ if and only if $f\in\loq$
and $\{R_{\ell}(f)\}_{\ell=1}^D\st\loq$.
Moreover, there exists a positive constant $C$
such that, for all $f\in\tls$,
$$
\frac1C\|f\|_{\tls}\le\sum_{\ell=0}^D\|R_{\ell}(f)\|_{\loq}
\le C\|f\|_{\tls},
$$
where $R_0:={\rm Id}$ denotes the \emph{identity operator}.
\end{theorem}

\begin{remark}\label{ra.l}
If $D=1$, Theorem \ref{ta.h} is just \cite[Theorem 1.3]{lly}.
\end{remark}

Fefferman-Stein decomposition says, for some function space $A$, there exists some space $B$ satisfying $B\nsubseteq A$ such that,
for $f\in A$, there exist $f_{l}\in B$ such that
$$f=\sum\limits^{D}_{l=0}R_{l} f_{l}.$$
The functions in B have better properties than those in A. But a function in $A$ has been written as a linear combination of a function in $B$ and the $n$ images of  functions in $B$ under correspondingly the $n$ Riesz transformations.
Such a result brings certain conveniences in PDE and in harmonic analysis.
The following theorems \ref{ta.x} and \ref{ta.i}  tell us that we have also Fefferman-Stein decomposition for ${\dot{F}^0_{\fz,\,q}(\rr^D)}$.

By Remark \ref{ra.k}(iv), we know that,
for any $q\in[2,\fz)$, $\tls\st\loq$ and ${\rm WE}^{\fz,\,q'}(\rr^D)\st\dtl$.
The following conclusions indicate that the above inclusions of sets are proper,
which are extensions of \cite[Remark 1.8]{lly}. The proof of theorem \ref{ta.x} will be given at section 3.

\begin{theorem}\label{ta.x}
Let $D\in\nn$ and $q\in[2,\fz)$. Then

(i) $\tls\subsetneqq L^{1}(\rr^D)  \bigcup \tls \subset {\rm WE}^{1,\,q}(\rr^D)\subset L^{1}(\rr^D) +    \tls$;

(ii) $ {\rm WE}^{\fz,\,q'}(\rr^D)\subsetneqq\dtl.$
\end{theorem}

Combining Theorem \ref{ta.h}, Remark \ref{ra.k}(iv)
and some arguments analogous to those
used in the proof of \cite[Theorem 1.7]{lly},
we obtain the following Fefferman-Stein decomposition
of ${\dot{F}^0_{\fz,\,q}(\rr^D)}$, the proof will be given in the final section.

\begin{theorem}\label{ta.i}
Let $D\in\nn$ and $q\in(1,2]$. Then $f\in{\dot{F}^0_{\fz,\,q}(\rr^D)}$
if and only if there exist $\{f_\ell\}_{\ell=0}^D\in \liq
$ such that
$f=f_0+\sum_{\ell=1}^D R_{\ell}\lf(f_\ell\r).$
\end{theorem}

By Theorems \ref{th:111} and \ref{ta.i}, we have
\begin{theorem}\label{ta.cor}
Let $D\in\nn$ and $q\in[2,\fz)$.
Then $f\in\schid$ belongs to $\tls$ if and only if $f\in L^{1}(\mathbb{R}^{D}) + \dot{F}^{0}_{1,q}(\mathbb{R}^{D})$
and $\{R_{\ell}(f)\}_{\ell=1}^D\st L^{1}(\mathbb{R}^{D}) + \dot{F}^{0}_{1,q}(\mathbb{R}^{D})$.
Moreover, there exists a positive constant $C$
such that, for all $f\in\tls$,
$$
\frac1C\|f\|_{\tls}\le\sum_{\ell=0}^D\|R_{\ell}(f)\|_{L^{1}(\mathbb{R}^{D}) + \dot{F}^{0}_{1,q}(\mathbb{R}^{D})}
\le C\|f\|_{\tls},
$$
where $R_0:={\rm Id}$ denotes the \emph{identity operator}.
\end{theorem}

\begin{remark}\label{ra.m}
Theorem \ref{ta.i} when $D=1$ is just \cite[Theorem 1.7]{lly}.
\end{remark}

The organization of this article is as follows.

In Section \ref{s2}, via the definition of the space ${\rm WE}^{1,\,q}(\rr^D)$,
the boundedness of Riesz transforms on $\tls$,
the Riesz transform characterization of $\hon$ and
some ideas from \cite{lly,yql},
we prove Theorem \ref{ta.h}, namely, establish the Riesz transform characterization
of Triebel-Lizorkin spaces $\dot{F}^0_{1,\,q}(\mathbb{R}^D)$.
Comparing with the corresponding proof of \cite[Subsection 6.2]{yql},
the main innovation of this proof is that we
regard the corresponding parts of the norms of Riesz transforms
$\{R_\ell(f_{s_1,\, N_1})\}_{\ell=1}^D$
in ${\rm WE}^{1,\,q}(\rr^D)$ as a whole to choose
$t^1_{s,\, N}\in\{0,\ldots, N+1\}$
such that \eqref{b.f} below holds true, while not to choose
$t^\ell_{s,\, N}\in\{0,\ldots, N+1\}$
such that \eqref{b.f} below holds true for each
$\ell\in\{1,\ldots,D\}$ separately as in \cite[(6.6)]{yql}.
Using this technique, we successfully overcome those difficulties described in
\cite[Remark 1.4]{lly}.

In Section \ref{s3}, we prove Theorem \ref{ta.x}. To this end, we first
give a $1$-dimensional Meyer wavelet satisfying $\psi(0)\neq0$
(see Example \ref{ec.a} below), which is taken from
\cite[Exercise 3.2]{w97}. By using such a $1$-dimensional Meyer
wavelet satisfying $\psi(0)\neq0$, we then finish the
proof of Theorem \ref{ta.x} via tensor products
and some arguments from the proof of \cite[Remark 1.8]{lly}.
Comparing with that proof of \cite[Remark 1.8]{lly}, we make
an additional assumption that $\psi(0)\neq0$ here, which is
needed in the estimate \eqref{3.4x} below.

In Section 4, we give  the proof of Theorems \ref{th:111}, \ref{ta.i} and \ref{ta.cor}.

Finally, we make some conventions on notation.
Throughout the whole paper, $C$ stands for a {\it positive constant} which
is independent of the main parameters, but it may vary from line to
line. If, for two real functions $f$ and $g$,
$f\le Cg$, we then write $f\ls g$;
if $f\ls g\ls f$, we then write $f\sim g$.
For $q\in(1,\fz)$, let $q'$ be the \emph{conjugate number} of $q$ defined
by $1/q+1/q'=1$. Let $\mathcal{C}$ be the set of complex numbers and $\nn:=\{1,2,\ldots\}$.
Furthermore, $\langle\cdot,\cdot\rangle$ and $(\cdot,\cdot)$
represent the duality relation, respectively, the $\ltw$ inner product.

\section{Proof of Theorem \ref{ta.h}}\label{s2}

\hskip\parindent In this section, we prove Theorem \ref{ta.h}.
To this end, we need to recall some well known results.

The following conclusion is taken from \cite[Corollary (8.21)]{fjw}.

\begin{theorem}\label{ta.f}
Let $D\in\nn$ and $q\in(1,\fz)$. Then the Riesz transform $R_{\ell}$
for each $\ell\in\{1,\ldots,D\}$ is bounded on $\tls$.
\end{theorem}

\begin{remark}\label{ra.p}
From Remark \ref{ra.n}(ii) and Theorem \ref{ta.f}, it follows that
the Riesz transform $R_{\ell}$
for each $\ell\in\{1,\ldots,D\}$ is bounded on $\hon$.
\end{remark}

The Riesz transform characterization of $\hon$ can be found
in \cite[p.\,221]{s70}.

\begin{theorem}\label{ta.g}
Let $D\in\nn$. The space $\hon$ is isomorphic to
the space of all functions $f\in\lon$
such that $\{R_{\ell}(f)\}_{\ell=1}^D\st\lon$.
Moreover, there exists a positive constant $C$ such that, for all $f\in\hon$,
$$
\frac1C\|f\|_{\hon}\le\|f\|_{\lon}+\sum_{\ell=1}^D\|R_{\ell}(f)\|_{\lon}
\le C\|f\|_{\hon}.
$$
\end{theorem}

The following lemma is completely analogous to \cite[Lemma 2.2]{lly},
the details being omitted.

\begin{lemma}\label{lb.j}
Let $D\in\nn$ and $q\in[2,\fz)$. If $f\in\loq$, then, for any $j\in\zz$,
$Q_{j}(f)\in\hon$, where
$Q_{j}(f):=\sum_{(\lz,\,k)\in E_D\times\zz^D}a^{\lz}_{j,\,k}(f)\psi^{\lz}_{j,\,k}$.
Moreover, there exists a positive constant $C$
such that, for all $j\in\zz$ and $f\in\loq$,
$$
\lf\|Q_j(f)\r\|_{\hon}\le C\|f\|_{\loq}.
$$
\end{lemma}

\begin{proof}[Proof of Theorem \ref{ta.h}]
We first show the necessity of Theorem \ref{ta.h}. By Remark \ref{ra.k}(iv)
and Theorem \ref{ta.f}, we have
$$
\sum_{\ell=0}^D\|R_\ell(f)\|_{\loq}\ls\sum_{\ell=0}^D\|R_\ell(f)\|_{\tls}
\ls\|f\|_{\tls},
$$
which completes the proof of the necessity of Theorem \ref{ta.h}.

Now we show the sufficiency of Theorem \ref{ta.h}. To this end,
for any $f\in\loq$ such that $\{R_{\ell}(f)\}_{\ell=1}^D\st\loq$,
it suffices to show that, for any
$s_1\in\zz$, $N_1\in\nn$ and $f_{s_1,\,N_1}:=P_{s_1,\,N_1}f$
defined as in \eqref{b.a}, we have
\begin{equation}\label{b.b}
\lf\|f_{s_1,\,N_1}\r\|_{\tls}\ls\sum_{\ell=0}^D
\lf\|R_{\ell}\lf(f_{s_1,\,N_1}\r)\r\|_{\loq},
\end{equation}
where the implicit constant is independent of $s_1$, $N_1$ and $f$.

Indeed, assume that \eqref{b.b} holds true for the time being.
Owing to \eqref{b.d}, for any $\ell\in\{1,\ldots,D\}$,
there exists a sequence $\{f^{\lz,\,\ell}_{j,\,k}\}_{(\lz,\,j,\,k)\in\blz_D}
\st\cc$ such that
$$
R_{\ell}\lf(f_{s_1,\,N_1}\r)
:=\sum_{\{(\lz,\,j,\,k)\in\blz_D:\ s_1-N_1-1\le j\le s_1+1\}}
f^{\lz,\,\ell}_{j,\,k}\psi^{\lz}_{j,\,k}\quad{\rm in}\quad \schid.
$$
By this and the orthogonality of $\{\psi^{\lz}_{j,\,k}\}_{(\lz,\,j,\,k)\in\blz_D}$,
we know that, for each $\ell\in\{1,\ldots,D\}$,
\begin{align}\label{b.x1}
&\lf\|R_{\ell}\lf(f_{s_1,\,N_1}\r)\r\|_{\loq}\\
&\noz\hs=\sup_{\{\wz{s}\in\zz,\,\wz{N}\in\nn\}}\min_{t\in\{0,\ldots,\wz{N}+1\}}
\lf[\lf\|T^{(1)}_{\wz{s},\,t,\,\wz{N}}R_{\ell}\lf(f_{s_1,\,N_1}\r)\r\|_{\tls}\r.\\
&\noz\hs\hs\lf.+\lf\|T^{(2)}_{\wz{s},\,t,\,\wz{N}}
R_{\ell}\lf(f_{s_1,\,N_1}\r)\r\|_{\lon}\r]\\
&\noz\hs=\sup_{\gfz{\wz{s}\in\zz,\,\wz{N}\in\nn}{\wz{s}\le s_1+1,\,
\wz{s}-\wz{N}\ge s_1-N_1-1}}
\min_{t\in\{0,\ldots,\wz{N}+1\}}
\lf[\lf\|T^{(1)}_{\wz{s},\,t,\,\wz{N}}R_{\ell}\lf(f_{s_1,\,N_1}\r)\r\|_{\tls}\r.\\
&\noz\hs\hs\lf.+\lf\|T^{(2)}_{\wz{s},\,t,\,\wz{N}}
R_{\ell}\lf(f_{s_1,\,N_1}\r)\r\|_{\lon}\r]\\
&\noz\hs=\sup_{\gfz{\wz{s}\in\zz,\,\wz{N}\in\nn}{\wz{s}\le s_1+1,\,
\wz{s}-\wz{N}\ge s_1-N_1-1}}
\min_{t\in\{0,\ldots,\wz{N}+1\}}
\lf[\lf\|T^{(1)}_{\wz{s},\,t,\,\wz{N}}R_{\ell}(f)\r\|_{\tls}
+\lf\|T^{(2)}_{\wz{s},\,t,\,\wz{N}}
R_{\ell}(f)\r\|_{\lon}\r]\\
&\noz\hs\le\lf\|R_{\ell}(f)\r\|_{\loq}<\fz
\end{align}
and, similarly,
\begin{align}\label{b.x2}
&\lf\|f_{s_1,\,N_1}\r\|_{\loq}\\
&\noz\hs=\sup_{\gfz{\wz{s}\in\zz,\,\wz{N}\in\nn}{\wz{s}\le s_1+1,\,
\wz{s}-\wz{N}\ge s_1-N_1-1}}
\min_{t\in\{0,\ldots,\wz{N}+1\}}
\lf[\lf\|T^{(1)}_{\wz{s},\,t,\,\wz{N}}\lf(f_{s_1,\,N_1}\r)\r\|_{\tls}\r.\\
&\noz\hs\hs\lf.+\lf\|T^{(2)}_{\wz{s},\,t,\,\wz{N}}
\lf(f_{s_1,\,N_1}\r)\r\|_{\lon}\r]\\
&\noz\hs=\sup_{\gfz{\wz{s}\in\zz,\,\wz{N}\in\nn}{\wz{s}\le s_1,\,
\wz{s}-\wz{N}\ge s_1-N_1}}
\min_{t\in\{0,\ldots,\wz{N}+1\}}
\lf[\lf\|T^{(1)}_{\wz{s},\,t,\,\wz{N}}(f_{s_1,\,N_1})\r\|_{\tls}
+\lf\|T^{(2)}_{\wz{s},\,t,\,\wz{N}}(f_{s_1,\,N_1})\r\|_{\lon}\r]\\
&\noz\hs=\sup_{\gfz{\wz{s}\in\zz,\,\wz{N}\in\nn}{\wz{s}\le s_1,\,
\wz{s}-\wz{N}\ge s_1-N_1}}
\min_{t\in\{0,\ldots,\wz{N}+1\}}
\lf[\lf\|T^{(1)}_{\wz{s},\,t,\,\wz{N}}(f)\r\|_{\tls}
+\lf\|T^{(2)}_{\wz{s},\,t,\,\wz{N}}(f)\r\|_{\lon}\r]\\
&\noz\hs\le\|f\|_{\loq}<\fz.
\end{align}
From \eqref{b.x1}, \eqref{b.x2} and \eqref{b.b}, we deduce that
$$
\lf\|f_{s_1,\,N_1}\r\|_{\tls}\ls\sum_{\ell=0}^D
\lf\|R_{\ell}(f)\r\|_{\loq}.
$$
This, together with Theorem \ref{ta.d} and the Levi lemma, implies that
$f\in\tls$ and
\begin{align*}
\|f\|_{\tls}
&\ls\lf\|\lf\{\sum_{(\lz,\,j,\,k)\in\blz_D}
\lf[2^{Dj}\lf|a^{\lz}_{j,\,k}(f)\r|\chi\lf(2^jx-k\r)\r]^q
\r\}^{1/q}\r\|_{\lon}\\
&\sim\lim_{N_1,s_1\to\fz}\lf\|\lf\{\sum_{\{(\lz,\,j,\,k)\in\blz_D:\
s_1-N_1\le j\le s_1\}}
\lf[2^{Dj}\lf|a^{\lz}_{j,\,k}(f)\r|\chi\lf(2^jx-k\r)\r]^q
\r\}^{1/q}\r\|_{\lon}\\
&\sim\lim_{N_1,\,s_1\to\fz}\lf\|f_{s_1,\,N_1}\r\|_{\tls}
\ls\sum_{\ell=0}^D\lf\|R_{\ell}(f)\r\|_{\loq},
\end{align*}
which are the desired conclusions.

Thus, to finish the proof of the sufficiency of Theorem \ref{ta.h},
we still need to prove \eqref{b.b}. To this end, fix $s_1\in\zz$ and $N_1\in\nn$.
In order to obtain the $\loq$-norms of $\{R_{\ell}(f_{s_1,\,N_1})\}_{\ell=0}^D$,
by \eqref{b.x1} and \eqref{b.x2}, it suffices to consider
$s:=s_1+1$ and $N:=N_1+2$ in \eqref{b.x} and \eqref{b.y}.
For such $s$ and $N$, there exist $t^{(0)}_{s,\,N},\,t^{(1)}_{s,\,N}
\in\{0,\ldots,N+1\}$ such that
\begin{align}\label{b.g}
&\lf\|T^{(1)}_{s,\,t^{(0)}_{s,\,N},\,N}\lf(f_{s_1,\,N_1}\r)\r\|_{\tls}
+\lf\|T^{(2)}_{s,\,t^{(0)}_{s,\,N},\,N}\lf(f_{s_1,\,N_1}\r)\r\|_{\lon}\\
&\noz\hs=\min_{t\in\{0,\ldots,N+1\}}
\lf[\lf\|T^{(1)}_{s,\,t,\,N}\lf(f_{s_1,\,N_1}\r)\r\|_{\tls}
+\lf\|T^{(2)}_{s,\,t,\,N}\lf(f_{s_1,\,N_1}\r)\r\|_{\lon}\r]
\end{align}
and
\begin{align}\label{b.f}
&\sum_{\ell=1}^D\lf[\lf\|T^{(1)}_{s,\,t^{(1)}_{s,\,N},\,N}R_{\ell}
\lf(f_{s_1,\,N_1}\r)\r\|_{\tls}
+\lf\|T^{(2)}_{s,\,t^{(1)}_{s,\,N},\,N}R_{\ell}\lf(f_{s_1,\,N_1}\r)\r\|_{\lon}\r]\\
&\noz\hs=\min_{t\in\{0,\ldots,N+1\}}\sum_{\ell=1}^D
\lf[\lf\|T^{(1)}_{s,\,t,\,N}R_{\ell}\lf(f_{s_1,\,N_1}\r)\r\|_{\tls}
+\lf\|T^{(2)}_{s,\,t,\,N}R_{\ell}\lf(f_{s_1,\,N_1}\r)\r\|_{\lon}\r].
\end{align}

In the remainder of this proof, to simplify the notation,
\emph{we let $g_1:=f_{s_1,\,N_1}$
for any fixed $s_1$ and $N_1$, $t_j:=t^{(j)}_{s,\,N}$ and
$T_{i,\,j}:=T^{(i)}_{s,\,t^{(j)}_{s,\,N},\,N}$ for any $i\in\{1,2\}$ and
$j\in\{0,1\}$}.

We then consider the following three cases.

\textbf{Case I}. $t_0=t_1$. In this case, we write
$g_1=a_1+a_2$, where
$$
a_1:=\sum_{j=s-t_0+1}^s Q_j\lf(g_1\r)
\quad {\rm and}\quad a_2:=\sum_{j=s-N}^{s-t_0}
Q_j\lf(g_1\r).
$$
By \eqref{b.g}, we have
$a_2=T_{2,\,0}(g_1)\in\lon$ and
\begin{equation}\label{x.o}
\lf\|a_2\r\|_{\lon}
=\lf\|T_{2,\,0}\lf(g_1\r)\r\|_{\lon}
\le\|g_1\|_{\loq},
\end{equation}
which, together with Lemma \ref{lb.j} and $\hon\st\lon$,
further implies that
$$
Q_{s-t_0}\lf(g_1\r)
+Q_{s-t_0-1}\lf(g_1\r)\in\hon
$$
and
\begin{align}\label{b.h}
&\lf\|Q_{s-t_0}\lf(g_1\r)
+Q_{s-t_0-1}\lf(g_1\r)\r\|_{\hon}\\
&\noz\hs\le \lf\|Q_{s-t_0}\lf(g_1\r)\r\|_{\hon}
+\lf\|Q_{s-t_0-1}\lf(g_1\r)\r\|_{\hon}
\ls\lf\|g_1\r\|_{\loq}.
\end{align}
Thus, by this, $\hon\st\lon$ and \eqref{x.o}, we obtain
$$
a_2-\lf[Q_{s-t_0}\lf(g_1\r)
+Q_{s-t_0-1}\lf(g_1\r)\r]\in\lon
$$
and
\begin{align}\label{b.i}
&\lf\|a_2-\lf[Q_{s-t_0}\lf(g_1\r)
+Q_{s-t_0-1}\lf(g_1\r)\r]\r\|_{\lon}\\
&\noz\hs\le\lf\|a_2\r\|_{\lon}
+\lf\|Q_{s-t_0}\lf(g_1\r)
+Q_{s-t_0-1}\lf(g_1\r)\r\|_{\hon}\ls\lf\|g_1\r\|_{\loq}.
\end{align}
Moreover, for each $\ell\in\{1,\ldots,D\}$, we have
\begin{align}\label{x.y}
T_{2,\,1}R_{\ell}\lf(g_1\r)
&=T_{2,\,1}R_{\ell}\lf(a_2
+Q_{s-t_0+1}\lf(g_1\r)\r)\\
&\noz=T_{2,\,1}R_{\ell}\lf(a_2
-\lf[Q_{s-t_0}\lf(g_1\r)+Q_{s-t_0-1}\lf(g_1\r)\r]\r)\\
&\noz\hs+T_{2,\,1}R_{\ell}\lf(Q_{s-t_0}\lf(g_1\r)
+Q_{s-t_0-1}\lf(g_1\r)\r)+T_{2,\,1}R_{\ell}Q_{s-t_0+1}\lf(g_1\r)\\
&\noz=R_{\ell}\lf(a_2
-\lf[Q_{s-t_0}\lf(g_1\r)+Q_{s-t_0-1}\lf(g_1\r)\r]\r)\\
&\noz\hs+T_{2,\,1}R_{\ell}\lf(Q_{s-t_0}\lf(g_1\r)
+Q_{s-t_0-1}\lf(g_1\r)\r)+T_{2,\,1}R_{\ell}Q_{s-t_0+1}\lf(g_1\r).
\end{align}
Hence, by \eqref{x.y}, \eqref{b.f}, $\hon\st\lon$,
Remarks \ref{ra.o} and \ref{ra.p}, \eqref{b.h} and Lemma \ref{lb.j},
we conclude that, for any $\ell\in\{1,\ldots,D\}$,
\begin{align*}
{\rm II}^{(\ell)}:&=R_{\ell}\lf(a_2
-\lf[Q_{s-t_0}\lf(g_1\r)+Q_{s-t_0-1}\lf(g_1\r)\r]\r)
+T_{2,\,1}R_{\ell}Q_{s-t_0+1}\lf(g_1\r)\in\lon
\end{align*}
and
\begin{align}\label{b.j}
\lf\|{\rm II}^{(\ell)}\r\|_{\lon}
&\le\lf\|T_{2,\,1}R_{\ell}\lf(g_1\r)\r\|_{\lon}\\
&\noz\hs+\lf\|T_{2,\,1}R_{\ell}
\lf(Q_{s-t_0}\lf(g_1\r)
+Q_{s-t_0-1}\lf(g_1\r)+Q_{s-t_0+1}\lf(g_1\r)\r)\r\|_{\hon}\\
&\noz\ls\sum_{\ell=1}^D\lf\|R_{\ell}\lf(g_1\r)\r\|_{\loq}\\
&\noz\hs+\lf\|R_{\ell}\lf(Q_{s-t_0}\lf(g_1\r)
+Q_{s-t_0-1}\lf(g_1\r)+Q_{s-t_0+1}\lf(g_1\r)\r)\r\|_{\hon}\\
&\noz\ls\sum_{\ell=1}^D\lf\|R_{\ell}\lf(g_1\r)\r\|_{\loq}\\
&\noz\hs+\lf\|Q_{s-t_0}\lf(g_1\r)
+Q_{s-t_0-1}\lf(g_1\r)+Q_{s-t_0+1}\lf(g_1\r)\r\|_{\hon}\\
&\noz\ls\sum_{\ell=1}^D\lf\|R_{\ell}\lf(g_1\r)\r\|_{\loq}+
\lf\|g_1\r\|_{\loq}.
\end{align}

From \eqref{b.d}, it follows that, for each $\ell\in\{1,\ldots,D\}$,
there exist $\{\tau^{\lz,\,\ell}_{j,\,k}\}_{(\lz,\,j,\,k)\in\blz_D}\st\cc$
such that
\begin{align*}
{\rm I}^{(\ell)}:&=R_{\ell}\lf(a_2
-\lf[Q_{s-t_0}\lf(g_1\r)
+Q_{s-t_0-1}\lf(g_1\r)\r]\r)
=\sum_{\{(\lz,\,j,\,k)\in\blz_D:\ s-N-1\le j\le s-t_0-1\}}
\tau^{\lz,\,\ell}_{j,\,k}\psi^{\lz}_{j,\,k}
\end{align*}
and
$$
R_{\ell}Q_{s-t_0+1}\lf(g_1\r)
=\sum_{\{(\lz,\,j,\,k)\in\blz_D:\ s-t_0\le j\le s-t_0+2\}}
\tau^{\lz,\,\ell}_{j,\,k}\psi^{\lz}_{j,\,k}.
$$

For any $h\in\li$ and $j_0\in\zz$, let
$$
P_{j_0}(h):=\sum_{k\in\zz^D}\lf\langle h,\psi^{\vec{0}}_{j_0,\,k}
\r\rangle\psi^{\vec{0}}_{j_0,\,k},
$$
where $\langle\cdot,\cdot\rangle$ represents the duality
between $\li$ and $\lon$.
We claim that $P_{j_0}(h)\in\li$. Indeed, by
$
|\langle h,\psi^{\vec{0}}_{j_0,\,k}\rangle|\ls2^{-Dj_0/2}
$
and $\psi^{\vec{0}}\in\sch$, we know that, for all $x\in\rr^D$,
$$
\lf|P_{j_0}(h)(x)\r|\ls\sum_{k\in\zz^D}2^{-Dj_0/2}
\lf|\psi^{\vec{0}}_{j_0,\,k}(x)\r|
\ls\sum_{k\in\zz^D}\lf|\psi^{\vec{0}}\lf(2^{j_0}x-k\r)\r|\ls1.
$$
Let
$$
h_0:=P_{s-t_0}(h)=\sum_{\{(\lz,\,j,\,k)\in\blz_D:\
j\le s-t_0-1\}}a^{\lz}_{j,\,k}(h)\psi^{\lz}_{j,\,k}.
$$
Thus, $h_0\in\li$ and $\|h_0\|_{\li}\ls1$ by the above claim.

Moreover, from \eqref{b.j}, we observe that, for any $\ell\in\{1,\ldots,D\}$,
$$
\lf|\lf\langle {\rm I}^{(\ell)},h\r\rangle\r|
=\lf|\lf\langle {\rm I}^{(\ell)},h_0\r\rangle\r|
=\lf|\lf\langle {\rm II}^{(\ell)},h_0\r\rangle\r|
\le\lf\|{\rm II}^{(\ell)}\r\|_{\lon}\lf\|h_0\r\|_{\li},
$$
which, combined with $\|h_0\|_{\li}\ls1$ and \eqref{b.j}, further implies that
$$
\lf\|{\rm I}^{(\ell)}\r\|_{\lon}\ls\lf\|{\rm II}^{(\ell)}\r\|_{\lon}
\ls\sum_{\ell=0}^D\lf\|R_{\ell}\lf(g_1\r)\r\|_{\loq}.
$$
From this, Theorem \ref{ta.g} and \eqref{b.i}, it follows that
$$
a_2-\lf[Q_{s-t_0}\lf(g_1\r)
+Q_{s-t_0-1}\lf(g_1\r)\r]\in\hon
$$
and
\begin{align*}
&\lf\|a_2-\lf[Q_{s-t_0}\lf(g_1\r)
+Q_{s-t_0-1}\lf(g_1\r)\r]\r\|_{\hon}\\
&\hs\sim\lf\|a_2-\lf[Q_{s-t_0}\lf(g_1\r)
+Q_{s-t_0-1}\lf(g_1\r)\r]\r\|_{\lon}\\
&\hs\hs+\sum_{\ell=0}^D\lf\|R_{\ell}\lf(a_2
-\lf[Q_{s-t_0}\lf(g_1\r)
+Q_{s-t_0-1}\lf(g_1\r)\r]\r)\r\|_{\lon}\\
&\hs\ls\lf\|g_1\r\|_{\loq}+
\sum_{\ell=1}^D\lf\|R_{\ell}\lf(g_1\r)\r\|_{\loq},
\end{align*}
which, together with Remark \ref{ra.n}(ii), \eqref{b.h}
and Lemma \ref{lb.j}, further implies that
\begin{align}\label{b.l}
\lf\|a_2\r\|_{\tls}&\ls\lf\|a_2\r\|_{\hon}\\
&\noz\ls\lf\|a_2-\lf[Q_{s-t_0}\lf(g_1\r)
+Q_{s-t_0-1}\lf(g_1\r)\r]\r\|_{\hon}\\
&\noz\hs+\lf\|Q_{s-t_0}\lf(g_1\r)
+Q_{s-t_0-1}\lf(g_1\r)\r\|_{\hon}\\
&\noz\ls\lf\|g_1\r\|_{\loq}+
\sum_{\ell=1}^D\lf\|R_{\ell}\lf(g_1\r)\r\|_{\loq}.
\end{align}

Furthermore, by \eqref{b.g}, we find that
$$
\lf\|a_1\r\|_{\tls}=\lf\|T_{1,\,0}
\lf(g_1\r)\r\|_{\tls}\le\lf\|g_1\r\|_{\tls},
$$
which, combined with \eqref{b.l}, implies that
$g_1=a_1+a_2\in\tls$ and
$$
\lf\|g_1\r\|_{\tls}\le\lf\|a_1\r\|_{\tls}
+\lf\|a_2\r\|_{\tls}
\ls\sum_{\ell=0}^D\lf\|R_{\ell}\lf(g_1\r)\r\|_{\loq}.
$$
This finishes the proof of \textbf{Case I}.

\textbf{Case II}. $t_0>t_1$. In this case, we write
$g_1=b_1+b_2+b_3$, where
$$
b_1:=\sum_{j=s-t_1+1}^s Q_j\lf(g_1\r),
\quad b_2:=\sum_{j=s-t_0+1}^{s-t_1}
Q_j\lf(g_1\r) \quad {\rm and}\quad b_3:=\sum_{j=s-N}^{s-t_0}
Q_j\lf(g_1\r).
$$
Similar to \eqref{x.y}, for any $\ell\in\{1,\ldots,D\}$, we know that
\begin{align*}
T_{2,\,1}R_{\ell}\lf(g_1\r)
&=T_{2,\,1}R_{\ell}\lf(b_3+b_2
+Q_{s-t_1+1}\lf(g_1\r)\r)\\
&=R_{\ell}\lf(b_3+b_2
-\lf[Q_{s-t_1}\lf(g_1\r)+Q_{s-t_1-1}\lf(g_1\r)\r]\r)\\
&\hs+T_{2,\,1}R_{\ell}\lf(Q_{s-t_1}\lf(g_1\r)
+Q_{s-t_1-1}\lf(g_1\r)\r)
+T_{2,\,1}R_{\ell}Q_{s-t_1+1}\lf(g_1\r)\\
&=:{\rm I}^{(\ell)}_1+ {\rm I}^{(\ell)}_2+{\rm I}^{(\ell)}_3.
\end{align*}

For any $h\in\li$, let
$
h_1:=\sum_{\{(\lz,\,j,\,k)\in\blz_D:\ j\le s-t_1-1\}}
a^{\lz}_{j,\,k}(h)\psi^{\lz}_{j,\,k}.
$
Similar to the proof of $h_0\in\li$, we have $h_1\in\li$ and
\begin{equation}\label{x.z}
\lf\|h_1\r\|_{\li}\ls1.
\end{equation}
By \eqref{b.d}, we know that, for any $\ell\in\{1,\ldots,D\}$,
there exists a sequence $\{f^{\lz,\,\ell}_{j,\,k}\}_{(\lz,\,j,\,k)\in\blz_D}
\st\cc$ such that
$${\rm I}^{(\ell)}_1=\sum_{\{(\lz,\,j,\,k)\in\blz_D:\ s-N-1\le j\le s-t_1-1\}}
f^{\lz,\,\ell}_{j,\,k}\psi^{\lz}_{j,\,k}\quad \mathrm{and}\quad
{\rm I}^{(\ell)}_3=\sum_{\{(\lz,\,j,\,k)\in\blz_D:\ j=s-t_1\}}
f^{\lz,\,\ell}_{j,\,k}\psi^{\lz}_{j,\,k},$$
which imply that
$$
\lf|\lf\langle {\rm I}^{(\ell)}_1,h\r\rangle\r|
=\lf|\lf\langle {\rm I}^{(\ell)}_1,h_1\r\rangle\r|
=\lf|\lf\langle {\rm I}^{(\ell)}_1+{\rm I}^{(\ell)}_3,h_1\r\rangle\r|
=\lf|\lf\langle T_{2,\,1}R_{\ell}\lf(g_1\r)
-{\rm I}^{(\ell)}_2,h_1\r\rangle\r|.
$$
Hence, by this, \eqref{x.z}, \eqref{b.f}, $\hon\st\lon$,
Remarks \ref{ra.o} and \ref{ra.p}, and Lemma \ref{lb.j},
we conclude that
\begin{align}\label{b.u}
\lf\|{\rm I}^{(\ell)}_1\r\|_{\lon}
&\ls\lf[\lf\|T_{2,\,1}R_{\ell}\lf(g_1\r)\r\|_{\lon}
+\lf\|{\rm I}^{(\ell)}_2\r\|_{\lon}\r]\\
&\noz\ls\sum_{\ell=1}^D\lf\|R_{\ell}\lf(g_1\r)\r\|_{\loq}
+\lf\|{\rm I}^{(\ell)}_2\r\|_{\hon}\\
&\noz\ls\sum_{\ell=1}^D\lf\|R_{\ell}\lf(g_1\r)\r\|_{\loq}
+\lf\|R_{\ell}\lf(Q_{s-t_1}\lf(g_1\r)
+Q_{s-t_1-1}\lf(g_1\r)\r)\r\|_{\hon}\\
&\noz\ls\sum_{\ell=1}^D\lf\|R_{\ell}\lf(g_1\r)\r\|_{\loq}
+\lf\|Q_{s-t_1}\lf(g_1\r)
+Q_{s-t_1-1}\lf(g_1\r)\r\|_{\hon}\\
&\noz\ls\sum_{\ell=1}^D\lf\|R_{\ell}\lf(g_1\r)\r\|_{\loq}+
\lf\|g_1\r\|_{\loq}.
\end{align}
Thus, ${\rm I}^{(\ell)}_1\in\lon$. Moreover, by Remark \ref{ra.p}
and Lemma \ref{lb.j}, we have
\begin{align}\label{x.u}
&\lf\|R_{\ell}\lf(Q_{s-t_1}\lf(g_1\r)
+Q_{s-t_1-1}\lf(g_1\r)\r)\r\|_{\hon}\\
&\noz\hs\ls\lf\|Q_{s-t_1}\lf(g_1\r)
+Q_{s-t_1-1}\lf(g_1\r)\r\|_{\hon}
\ls\lf\|g_1\r\|_{\loq}.
\end{align}
From this, $\hon\st\lon$ and \eqref{b.u}, we deduce that
\begin{align}\label{x.v}
\lf\|R_{\ell}\lf(b_2+b_3\r)\r\|_{\lon}
&\le\lf\|{\rm I}^{(\ell)}_1\r\|_{\lon}
+\lf\|R_{\ell}\lf(Q_{s-t_1}\lf(g_1\r)
+Q_{s-t_1-1}\lf(g_1\r)\r)\r\|_{\hon}\\
&\noz\ls\sum_{\ell=0}^D\lf\|R_{\ell}\lf(g_1\r)\r\|_{\loq}.
\end{align}

On the other hand, for any $h\in\li$, let
$
\wz{h}_0:=\sum_{\{(\lz,\,j,\,k)\in\blz_D:\ j\ge s-t_0+2\}}
a^{\lz}_{j,\,k}(h)\psi^{\lz}_{j,\,k}.
$
Similar to the proof of $h_0\in\li$, we have $h-\wz{h}_0\in\li$ and
$\|h-\wz{h}_0\|_{\li}\ls1$, which further implies that
\begin{equation}\label{x.w}
\lf\|\wz{h}_0\r\|_{\li}\le\lf\|h\r\|_{\li}
+\lf\|h-\wz{h}_0\r\|_{\li}\ls1.
\end{equation}
By \eqref{b.d}, we know that
\begin{align*}
&R_{\ell}\lf(b_2-\lf[Q_{s-t_0+1}\lf(g_1\r)
+Q_{s-t_0+2}\lf(g_1\r)\r]\r)
=\sum_{\{(\lz,\,j,\,k)\in\blz_D:\ s-t_0+2\le j\le s-t_1+1\}}
f^{\lz,\,\ell}_{j,\,k}\psi^{\lz}_{j,\,k},
\end{align*}
which implies that
\begin{align}\label{x.n}
&\lf\langle R_{\ell}\lf(b_2-\lf[Q_{s-t_0+1}\lf(g_1\r)
+Q_{s-t_0+2}\lf(g_1\r)\r]\r),h\r\rangle\\
&\noz\hs=\lf\langle R_{\ell}\lf(b_2-\lf[Q_{s-t_0+1}\lf(g_1\r)
+Q_{s-t_0+2}\lf(g_1\r)\r]\r),\wz{h}_0\r\rangle\\
&\noz\hs=\lf\langle R_{\ell}\lf(b_3+b_2
-\lf[Q_{s-t_0+1}\lf(g_1\r)
+Q_{s-t_0+2}\lf(g_1\r)\r]\r),\wz{h}_0\r\rangle.
\end{align}
From an argument similar to that used in \eqref{x.u}, it follows that
\begin{align}\label{x.t}
&\lf\|R_{\ell}\lf(Q_{s-t_0+1}\lf(g_1\r)
+Q_{s-t_0+2}\lf(g_1\r)\r)\r\|_{\hon}\\
&\noz\hs\ls\lf\|Q_{s-t_0+1}\lf(g_1\r)
+Q_{s-t_0+2}\lf(g_1\r)\r\|_{\hon}
\ls\lf\|g_1\r\|_{\loq}.
\end{align}
Thus, by \eqref{x.n}, \eqref{x.w}, \eqref{x.t},
$\hon\st\lon$ and \eqref{x.v}, we conclude that
\begin{align*}
&\lf\|R_{\ell}\lf(b_2-\lf[Q_{s-t_0+1}\lf(g_1\r)
+Q_{s-t_0+2}\lf(g_1\r)\r]\r)\r\|_{\lon}\\
&\hs\ls\lf\|R_{\ell}\lf(b_3+b_2
-\lf[Q_{s-t_0+1}\lf(g_1\r)
+Q_{s-t_0+2}\lf(g_1\r)\r]\r)\r\|_{\lon}\\
&\hs\ls\lf\|R_{\ell}\lf(b_3+b_2\r)\r\|_{\lon}
+\lf\|R_{\ell}\lf(Q_{s-t_0+1}\lf(g_1\r)
+Q_{s-t_0+2}\lf(g_1\r)\r)\r\|_{\hon}\\
&\hs\ls\sum_{\ell=0}^D\lf\|R_{\ell}\lf(g_1\r)\r\|_{\loq}.
\end{align*}
Therefore, by this, $\hon\st\lon$ and \eqref{x.t}, we obtain
\begin{align*}
\lf\|R_{\ell}\lf(b_2\r)\r\|_{\lon}
&\le\lf\|R_{\ell}\lf(b_2
-\lf[Q_{s-t_0+1}\lf(g_1\r)
+Q_{s-t_0+2}\lf(g_1\r)\r]\r)\r\|_{\lon}\\
&\hs+\lf\|R_{\ell}\lf(Q_{s-t_0+1}\lf(g_1\r)
+Q_{s-t_0+2}\lf(g_1\r)\r)\r\|_{\hon}\\
&\ls\sum_{\ell=0}^D\lf\|R_{\ell}\lf(g_1\r)\r\|_{\loq},
\end{align*}
which, together with \eqref{x.v}, implies that
 \begin{align}\label{x.s}
\lf\|R_{\ell}\lf(b_3\r)\r\|_{\lon}
&\le\lf\|R_{\ell}\lf(b_3+b_2\r)\r\|_{\lon}
+\lf\|R_{\ell}\lf(b_2\r)\r\|_{\lon}\\
&\noz\ls\sum_{\ell=0}^D\lf\|R_{\ell}\lf(g_1\r)\r\|_{\loq}.
\end{align}

Furthermore, by \eqref{b.g}, we know that
$$
\lf\|b_3\r\|_{\lon}=
\lf\|T_{2,\,0}\lf(g_1\r)\r\|_{\lon}
\ls\lf\|g_1\r\|_{\loq},
$$
which, combined with \eqref{x.s} and Theorem \ref{ta.g}, implies that
$b_3\in\hon$ and
\begin{align}\label{b.o}
\lf\|b_3\r\|_{\hon}
&\sim\lf\|b_3\r\|_{\lon}
+\sum_{\ell=1}^D\lf\|R_{\ell}\lf(b_3\r)\r\|_{\lon}
\ls\sum_{\ell=0}^D\lf\|R_{\ell}\lf(g_1\r)\r\|_{\loq}.
\end{align}
By \eqref{b.o} and Remark \ref{ra.n}(ii), we know that
$b_3\in\tls$ and
\begin{equation}\label{b.q}
\lf\|b_3\r\|_{\tls}\ls\lf\|b_3\r\|_{\hon}
\ls\sum_{\ell=0}^D\lf\|R_{\ell}\lf(g_1\r)\r\|_{\loq}.
\end{equation}

Moreover, by \eqref{b.g}, we obtain
$$
\lf\|b_1+b_2\r\|_{\tls}=
\lf\|T_{1,\,0}\lf(g_1\r)\r\|_{\tls}
\ls\lf\|g_1\r\|_{\loq},
$$
which, together with \eqref{b.q} and Remark \ref{ra.n}(ii),
further implies that
$$
\lf\|g_1\r\|_{\tls}\le\lf\|b_1+b_2\r\|_{\tls}
+\lf\|b_3\r\|_{\tls}
\ls\sum_{\ell=0}^D\lf\|R_{\ell}\lf(g_1\r)\r\|_{\loq}.
$$
This finishes the proof of \textbf{Case II}.

\textbf{Case III}. $t_0<t_1$. In this case, we write
$g_1=e_1+e_2+e_3$, where
$$
e_1:=\sum_{j=s-t_0+1}^s Q_j\lf(g_1\r),
\quad e_2:=\sum_{j=s-t_1+1}^{s-t_0}
Q_j\lf(g_1\r) \quad{\rm and}\quad e_3:=\sum_{j=s-N}^{s-t_1}
Q_j\lf(g_1\r).
$$
Similar to \eqref{x.y}, for any $\ell\in\{1,\ldots,D\}$, we have
\begin{align}\label{x.m}
T_{2,\,1}R_{\ell}\lf(g_1\r)
&=T_{2,\,1}R_{\ell}\lf(e_3
+Q_{s-t_1+1}\lf(g_1\r)\r)\\
&\noz=R_{\ell}\lf(e_3
-\lf[Q_{s-t_1}\lf(g_1\r)+Q_{s-t_1-1}\lf(g_1\r)\r]\r)\\
&\noz\hs+T_{2,\,1}R_{\ell}\lf(Q_{s-t_1}\lf(g_1\r)
+Q_{s-t_1-1}\lf(g_1\r)\r)+T_{2,\,1}R_{\ell}Q_{s-t_1+1}\lf(g_1\r)\\
&\noz=:{\rm II}^{(\ell)}_1+ {\rm II}^{(\ell)}_2+{\rm II}^{(\ell)}_3.
\end{align}
For any $h\in\li$, let
$$
h_2:=\sum_{\{(\lz,\,j,\,k)\in\blz_D:\ j\le s-t_1-1\}}
a^{\lz}_{j,\,k}(h)\psi^{\lz}_{j,\,k}.
$$
By an argument similar to that used in the proof of $h_0\in\li$,
we conclude that $h_2\in\li$ and
$\|h_2\|_{\li}\ls1$. By \eqref{b.d},
we know that, for any $\ell\in\{1,\ldots,D\}$,
there exists a sequence $\{f^{\lz,\,\ell}_{j,\,k}\}_{(\lz,\,j,\,k)\in\blz_D}
\st\cc$ such that
$$
{\rm II}^{(\ell)}_1=\sum_{\{(\lz,\,j,\,k)
\in\blz_D:\ s-N-1\le j\le s-t_1-1\}}
f^{\lz,\,\ell}_{j,\,k}\psi^{\lz}_{j,\,k}
$$
and
$$
T_{2,\,1}R_{\ell}Q_{s-t_1+1}\lf(g_1\r)
=\sum_{\{(\lz,\,j,\,k)\in\blz_D:\ j=s-t_1\}}
f^{\lz,\,\ell}_{j,\,k}\psi^{\lz}_{j,\,k},
$$
which, together with \eqref{x.m}, imply that
$$
\lf|\lf\langle {\rm II}^{(\ell)}_1,h\r\rangle\r|
=\lf|\lf\langle {\rm II}^{(\ell)}_1,h_2\r\rangle\r|
=\lf|\lf\langle {\rm II}^{(\ell)}_1+{\rm II}^{(\ell)}_3,h_2\r\rangle\r|
=\lf|\lf\langle T_{2,\,1}R_{\ell}\lf(g_1\r)
-{\rm II}^{(\ell)}_2,h_2\r\rangle\r|.
$$
Hence, by this, $\|h_2\|_{\li}\ls1$, \eqref{b.f}, $\hon\st\lon$,
Remarks \ref{ra.o} and \ref{ra.p}, and Lemma \ref{lb.j},
we conclude that
\begin{align}\label{b.m}
\lf\|{\rm II}^{(\ell)}_1\r\|_{\lon}
&\ls\lf\|T_{2,\,1}R_{\ell}\lf(g_1\r)\r\|_{\lon}
+\lf\|{\rm II}^{(\ell)}_2\r\|_{\lon}\\
&\noz\ls\sum_{\ell=1}^D\lf\|R_{\ell}\lf(g_1\r)\r\|_{\loq}
+\lf\|{\rm II}^{(\ell)}_2\r\|_{\hon}\\
&\noz\ls\sum_{\ell=1}^D\lf\|R_{\ell}\lf(g_1\r)\r\|_{\loq}
+\lf\|R_{\ell}\lf(Q_{s-t_1}\lf(g_1\r)
+Q_{s-t_1-1}\lf(g_1\r)\r)\r\|_{\hon}\\
&\noz\ls\sum_{\ell=1}^D\lf\|R_{\ell}\lf(g_1\r)\r\|_{\loq}
+\lf\|Q_{s-t_1}\lf(g_1\r)
+Q_{s-t_1-1}\lf(g_1\r)\r\|_{\hon}\\
&\noz\ls\sum_{\ell=1}^D\lf\|R_{\ell}\lf(g_1\r)\r\|_{\loq}+
\lf\|g_1\r\|_{\loq}.
\end{align}
Thus, ${\rm II}^{(\ell)}_1\in\lon$.

By \eqref{b.f}, we have $e_2+e_3\in\lon$ and
\begin{equation}\label{b.s}
\lf\|e_2+e_3\r\|_{\lon}
=\lf\|T_{2,\,0}\lf(g_1\r)\r\|_{\lon}
\ls\lf\|g_1\r\|_{\loq}.
\end{equation}
For any $h\in\li$, let
$$
\wz{h}_1:=\sum_{\{(\lz,\,j,\,k)\in\blz_D:\
j\le s-t_1-2\}}h^{\lz}_{j,\,k}\psi^{\lz}_{j,\,k}.
$$
Similar to the proof of $h_0\in\li$, we have $\|\wz{h}_1\|_{\li}\ls1$.
We notice that
\begin{align*}
&\lf\langle e_3-\lf[Q_{s-t_1}\lf(g_1\r)
+Q_{s-t_1-1}\lf(g_1\r)\r],h\r\rangle\\
&\hs=\lf\langle e_3-\lf[Q_{s-t_1}\lf(g_1\r)
+Q_{s-t_1-1}\lf(g_1\r)\r],\wz{h}_1\r\rangle
=\lf\langle e_3,\wz{h}_1\r\rangle
=\lf\langle e_3+e_2,\wz{h}_1\r\rangle.
\end{align*}
Therefore, by this, \eqref{b.s} and $\|\wz{h}_1\|_{\li}\ls1$, we obtain
\begin{align*}
\lf\|e_3-\lf[Q_{s-t_1}\lf(g_1\r)
+Q_{s-t_1-1}\lf(g_1\r)\r]\r\|_{\lon}&\ls\lf\|e_3+e_2\r\|_{\lon}
\ls\lf\|g_1\r\|_{\loq}.
\end{align*}
Hence $e_3-[Q_{s-t_1}(g_1)
+Q_{s-t_1-1}(g_1)]\in\lon$.
From this, \eqref{b.m} and Theorem \ref{ta.g}, we deduce that
$e_3-[Q_{s-t_1}(g_1)
+Q_{s-t_1-1}(g_1)]\in\hon$ and
\begin{align*}
&\lf\|e_3-\lf[Q_{s-t_1}\lf(g_1\r)
+Q_{s-t_1-1}\lf(g_1\r)\r]\r\|_{\hon}\\
&\hs\sim\lf\|e_3-\lf[Q_{s-t_1}\lf(g_1\r)
+Q_{s-t_1-1}\lf(g_1\r)\r]\r\|_{\lon}\\
&\hs\hs+\sum_{\ell=1}^D\lf\|R_{\ell}
\lf(e_3-\lf[Q_{s-t_1}\lf(g_1\r)
+Q_{s-t_1-1}\lf(g_1\r)\r]\r)\r\|_{\lon}\\
&\hs\ls\sum_{\ell=0}^D\lf\|R_{\ell}\lf(g_1\r)\r\|_{\loq}.
\end{align*}
Then, by this, $\hon\st\lon$ and Lemma \ref{lb.j},
we know that $e_3\in\lon$ and
\begin{align}\label{b.z}
\lf\|e_3\r\|_{\lon}
&\le\lf\|e_3-\lf[Q_{s-t_1}\lf(g_1\r)
+Q_{s-t_1-1}\lf(g_1\r)\r]\r\|_{\lon}\\
&\noz\hs+\lf\|Q_{s-t_1}\lf(g_1\r)
+Q_{s-t_1-1}\lf(g_1\r)\r\|_{\hon}\\
&\ls\sum_{\ell=0}^D\lf\|R_{\ell}\lf(g_1\r)\r\|_{\loq}.\noz
\end{align}
For each $\ell\in\{1,\ldots, D\}$, we observe that
\begin{align*}
T_{1,\,1}R_{\ell}\lf(g_1\r)
=T_{1,\,1}R_{\ell}\lf(e_1+e_2
+Q_{s-t_1}\lf(g_1\r)\r).
\end{align*}
By this and \eqref{b.f}, we know that, for any $\ell\in\{1,\ldots, D\}$,
$$
T_{1,\,1}R_{\ell}\lf(e_1+e_2
+Q_{s-t_1}\lf(g_1\r)\r)\in\tls
$$
and
$$
\lf\|T_{1,\,1}R_{\ell}\lf(e_1+e_2
+Q_{s-t_1}\lf(g_1\r)\r)\r\|_{\tls}\ls
\sum_{\ell=1}^D\lf\|R_{\ell}\lf(g_1\r)\r\|_{\loq}.
$$
This, together with
\begin{equation}\label{b.tx}
\lf\|e_1\r\|_{\tls}
=\lf\|T_{1,\,0}\lf(g_1\r)\r\|_{\tls}
\le\lf\|g_1\r\|_{\loq}\quad({\rm see}\ (\ref{b.g})),
\end{equation}
Theorems \ref{ta.d} and \ref{ta.f},
and \eqref{b.tx}, further implies that, for each $\ell\in\{1,\ldots,D\}$,
\begin{align}\label{b.v}
&\lf\|T_{1,\,1}R_{\ell}\lf(e_2
+Q_{s-t_0}\lf(g_1\r)\r)\r\|_{\tls}\\
&\noz\hs\le\lf\|T_{1,\,1}R_{\ell}\lf(e_1+e_2
+Q_{s-t_0}\lf(g_1\r)\r)\r\|_{\tls}
+\lf\|T_{1,\,1}R_{\ell}\lf(e_1\r)\r\|_{\tls}\\
&\noz\hs\ls\sum_{\ell=1}^D\lf\|R_{\ell}\lf(g_1\r)\r\|_{\loq}
+\lf\|R_{\ell}\lf(e_1\r)\r\|_{\tls}\\
&\noz\hs\ls\sum_{\ell=1}^D\lf\|R_{\ell}\lf(g_1\r)\r\|_{\loq}
+\lf\|e_1\r\|_{\tls}\\
&\noz\hs\ls\sum_{\ell=1}^D\lf\|R_{\ell}\lf(g_1\r)\r\|_{\loq}
+\lf\|g_1\r\|_{\loq}.
\end{align}

Furthermore, for any $\ell\in\{1,\ldots,D\}$, we notice that
\begin{align}\label{x.l}
T_{1,\,1}R_{\ell}\lf(e_2
+Q_{s-t_1}\lf(g_1\r)\r)&=R_{\ell}\lf(e_2
-\lf[Q_{s-t_1+1}\lf(g_1\r)+Q_{s-t_0}\lf(g_1\r)\r]\r)\\
&\noz\hs+T_{1,\,1}R_{\ell}
\lf(Q_{s-t_1+1}\lf(g_1\r)
+Q_{s-t_1}\lf(g_1\r)+Q_{s-t_0}\lf(g_1\r)\r).
\end{align}
By Theorems \ref{ta.d} and \ref{ta.f}, Remark \ref{ra.n}(ii)
and Lemma \ref{lb.j}, we conclude that
\begin{align*}
&\lf\|T_{1,\,1}R_{\ell}
\lf(Q_{s-t_1+1}\lf(g_1\r)
+Q_{s-t_1}\lf(g_1\r)+Q_{s-t_0}\lf(g_1\r)\r)\r\|_{\tls}\\
&\hs\ls\lf\|R_{\ell}\lf(Q_{s-t_1+1}\lf(g_1\r)
+Q_{s-t_1}\lf(g_1\r)+Q_{s-t_0}\lf(g_1\r)\r)\r\|_{\tls}\\
&\hs\ls\lf\|Q_{s-t_1+1}\lf(g_1\r)
+Q_{s-t_1}\lf(g_1\r)+Q_{s-t_0}\lf(g_1\r)\r\|_{\tls}\\
&\hs\ls\lf\|Q_{s-t_1+1}\lf(g_1\r)
+Q_{s-t_1}\lf(g_1\r)+Q_{s-t_0}\lf(g_1\r)\r\|_{\hon}
\ls\lf\|g_1\r\|_{\tls},
\end{align*}
which, together with \eqref{x.l} and \eqref{b.v}, implies that
\begin{align}\label{y.u}
&\lf\|R_{\ell}\lf(e_2
-\lf[Q_{s-t_1+1}\lf(g_1\r)
+Q_{s-t_0}\lf(g_1\r)\r]\r)\r\|_{\tls}\\
&\noz\hs\le\lf\|T_{1,\,1}R_{\ell}\lf(e_2
+Q_{s-t_1}\lf(g_1\r)\r)\r\|_{\tls}\\
&\noz\hs\hs+\lf\|T_{1,\,1}R_{\ell}
\lf(Q_{s-t_1+1}\lf(g_1\r)
+Q_{s-t_1}\lf(g_1\r)
+Q_{s-t_0}\lf(g_1\r)\r)\r\|_{\tls}\\
&\noz\hs\ls\lf\|g_1\r\|_{\loq}+
\sum_{\ell=1}^D\lf\|R_{\ell}\lf(g_1\r)\r\|_{\loq}.
\end{align}

Now we need a useful identity from \cite[p.\,224, (2.9)]{sw} that,
for all $f\in\ltw$,
\begin{equation}\label{y.x}
\sum_{\ell=1}^D R^2_\ell(f)=-f.
\end{equation}

From $e_2-[Q_{s-t_1+1}(g_1)
+Q_{s-t_0}(g_1)]\in\ltw$ and \eqref{y.x}, we deduce that
\begin{align*}
e_2-\lf[Q_{s-t_1+1}\lf(g_1\r)
+Q_{s-t_0}\lf(g_1\r)\r]=\sum_{\ell=1}^D R^2_{\ell}\lf(e_2
-\lf[Q_{s-t_1+1}\lf(g_1\r)+Q_{s-t_0}\lf(g_1\r)\r]\r)
\in\tls,
\end{align*}
which, combined with Theorem \ref{ta.f} and \eqref{y.u}, implies that
\begin{align}\label{y.y}
&\lf\|e_2-\lf[Q_{s-t_1+1}\lf(g_1\r)
+Q_{s-t_0}\lf(g_1\r)\r]\r\|_{\tls}\\
&\noz\hs\le\sum_{\ell=1}^D \lf\|R^2_{\ell}\lf(e_2
-\lf[Q_{s-t_1+1}\lf(g_1\r)
+Q_{s-t_0}\lf(g_1\r)\r]\r)\r\|_{\tls}\\
&\noz\hs\ls\sum_{\ell=1}^D \lf\|R_{\ell}\lf(e_2
-\lf[Q_{s-t_1+1}\lf(g_1\r)
+Q_{s-t_0}\lf(g_1\r)\r]\r)\r\|_{\tls}\ls\sum_{\ell=0}^D\lf\|R_{\ell}\lf(g_1\r)\r\|_{\loq}.
\end{align}

Again, by Remark \ref{ra.n}(ii) and Lemma \ref{lb.j}, we obtain
\begin{align*}
\lf\|Q_{s-t_1+1}\lf(g_1\r)
+Q_{s-t_0}\lf(g_1\r)\r\|_{\tls}\!\ls\lf\|Q_{s-t_1+1}\lf(g_1\r)
+Q_{s-t_0}\lf(g_1\r)\r\|_{\hon}
\ls\lf\|g_1\r\|_{\loq},
\end{align*}
which, together with \eqref{y.y}, implies that
\begin{align}\label{b.w}
\lf\|e_2\r\|_{\tls}
&\le\lf\|e_2-\lf[Q_{s-t_1+1}\lf(g_1\r)
+Q_{s-t_0}\lf(g_1\r)\r]\r\|_{\tls}\\
&\noz\hs+\lf\|Q_{s-t_1+1}\lf(g_1\r)
+Q_{s-t_0}\lf(g_1\r)\r\|_{\tls}
\ls\sum_{\ell=0}^D\lf\|R_{\ell}\lf(g_1\r)\r\|_{\loq}.
\end{align}

Combining with \eqref{b.tx}, \eqref{b.w} and \eqref{b.z}, we obtain
\begin{align*}
\lf\|g_1\r\|_{\tls}&\le\sum_{j=1}^3\lf\|e_j\r\|_{\tls}
\ls\sum_{\ell=0}^D\lf\|R_{\ell}\lf(g_1\r)\r\|_{\loq},
\end{align*}
which completes the proof of \textbf{Case III} and hence
Theorem \ref{ta.h}.
\end{proof}

\section{Proof of Theorem \ref{ta.x}}\label{s3}

\hskip\parindent In this section, we prove Theorem \ref{ta.x}.
To this end, we first provide a
$1$-dimensional Meyer wavelet satisfying $\psi(0)\neq0$,
which is taken from \cite[Exercise 3.2]{w97}.

\begin{example}\label{ec.a}
Let
$$
f(x):=\begin{cases}
e^{-1/x^2}, \ \  \ \ &x\in(0,\fz),\\
0, \ \  \ \ &x\in(-\fz,0],
\end{cases}
$$
$f_1(x):=f(x)f(1-x)$ and $g(x):=\lf[\int_{-\fz}^\fz f_1(t)\,dt\r]^{-1}
\int_{-\fz}^x f_1(t)\,dt$ for all $x\in\rr$.
Let $\xi\in\rr$ and $\Phi(\xi):=\frac1{\sqrt{2\pi}}\cos
(\frac{\pi}2 g(\frac{3}{2\pi}|\xi|-1))$. Then, by
\cite[Exercise 3.2]{w97}, we know that $\Phi\in C^\fz(\rr)$ satisfies
\eqref{ph1} through \eqref{ph5}.
Following the construction of the $1$-dimensional Meyer wavelet, we obtain
the ``father" wavelet $\phi$, the corresponding function $m_{\phi}$
of $\phi$, and the ``mother" wavelet $\psi$.

By the proof of
\cite[Proposition 3.3(ii)]{w97}, we know that, for any $x\in\rr$,
$$
\psi(x)=\frac1{\sqrt{2\pi}}\int_{-\fz}^\fz
\cos\lf([x+1/2]\xi\r)\az(\xi)\,d\xi,
$$
where, for each $\xi\in\rr$, $\az(\xi):=m_{\phi}(\xi/2+\pi)\Phi(\xi/2)$
is an even function supported in $[-8\pi/3,-2\pi/3]\cup[2\pi/3,8\pi/3]$.

Now we show that $\psi(0)\neq0$. Indeed, by the facts
that $\az$ is even, $\supp(\az)\st[-8\pi/3,-2\pi/3]\cup[2\pi/3,8\pi/3]$
and $m_{\phi}$ is $2\pi$-periodic, we have
\begin{align*}
\psi(0)&=\frac2{\sqrt{2\pi}}\int_0^\fz \cos(\xi/2)\az(\xi)\,d\xi
=\frac2{\sqrt{2\pi}}\int_{2\pi/3}^{8\pi/3}\cos(\xi/2)\az(\xi)\,d\xi\\
&=\frac2{\sqrt{2\pi}}\int_{2\pi/3}^{8\pi/3}\cos(\xi/2)
m_{\phi}(\xi/2+\pi-2\pi)\Phi(\xi/2)\,d\xi\\
&=2\int_{2\pi/3}^{8\pi/3}\cos(\xi/2)
\Phi(\xi-2\pi)\Phi(\xi/2)\,d\xi\\
&=2\int_{2\pi/3}^{\pi}\cos(\xi/2)\Phi(\xi-2\pi)\Phi(\xi/2)\,d\xi
+2\int_{\pi}^{4\pi/3}\cdots+2\int_{4\pi/3}^{8\pi/3}\cdots
=:{\rm J}_1+{\rm J}_2+{\rm J}_3.
\end{align*}

We first estimate ${\rm J}_1$. Observe that, for any $\xi\in[2\pi/3,\pi]$,
by \eqref{ph3}, we know that $\Phi(\xi/2)=1/\sqrt{2\pi}$.
Moreover, from $\frac12\le \frac{3}{2\pi}|\xi-2\pi|-1\le1$
and the construction of $g$, we deduce that
$g(\frac12)\le g(\frac{3}{2\pi}|\xi-2\pi|-1)\le g(1)=1$.
Hence, by the construction of $\Phi$, we obtain
$$
0\le\Phi(\xi-2\pi)=\frac1{2\pi}\cos\lf(\frac{\pi}2 g\lf(
\frac{3}{2\pi}|\xi-2\pi|-1\r)\r)\le\frac1{\sqrt{2\pi}}\cos\lf(\frac{\pi}2
g\lf(\frac12\r)\r),
$$
which further implies that
\begin{align*}
\lf|{\rm J}_1\r|&\le2\int_{2\pi/3}^{\pi}\cos(\xi/2)\frac1{\sqrt{2\pi}}
\cos\lf(\frac{\pi}2g\lf(\frac12\r)\r)\frac1{\sqrt{2\pi}}\,d\xi\\
&=\frac1{\pi}\cos\lf(\frac{\pi}2g\lf(\frac12\r)\r)
\int_{2\pi/3}^{\pi}\cos(\xi/2)\,d\xi
=\frac{2-\sqrt{3}}{\pi}\cos\lf(\frac{\pi}2g\lf(\frac12\r)\r).
\end{align*}

Now we deal with ${\rm J}_2$. Observe that, for any $\xi\in[\pi,4\pi/3]$,
by \eqref{ph3}, we know that $\Phi(\xi/2)=1/\sqrt{2\pi}$.
Moreover, from $0\le \frac{3}{2\pi}|\xi-2\pi|-1\le\frac12$
and the construction of $g$, it follows
that $0=g(0)\le g(\frac{3}{2\pi}|\xi-2\pi|-1)\le g(\frac12)$.
Hence, by the construction of $\Phi$ again, we have
$$
1\ge\Phi(\xi-2\pi)\ge\frac1{\sqrt{2\pi}}\cos\lf(\frac{\pi}2
g\lf(\frac12\r)\r),
$$
which further implies that
\begin{align*}
\lf|{\rm J}_2\r|&=-2\int^{4\pi/3}_{\pi}\cos(\xi/2)
\Phi(\xi-2\pi)\Phi(\xi/2)\,d\xi\\
&\ge-2\int^{4\pi/3}_{\pi}\cos(\xi/2)\frac1{\sqrt{2\pi}}
\cos\lf(\frac{\pi}2g\lf(\frac12\r)\r)\frac1{\sqrt{2\pi}}\,d\xi\\
&=-\frac1{\pi}\cos\lf(\frac{\pi}2g\lf(\frac12\r)\r)
\int^{4\pi/3}_{\pi}\cos(\xi/2)\,d\xi
=\frac{2-\sqrt{3}}{\pi}\cos\lf(\frac{\pi}2g\lf(\frac12\r)\r).
\end{align*}

Finally, we estimate ${\rm J}_3$. We first write
\begin{equation}\label{x.r}
\lf|{\rm J}_3\r|=-2\int^{8\pi/3}_{4\pi/3}\cos(\xi/2)
\Phi(\xi-2\pi)\Phi(\xi/2)\,d\xi
\ge-2\int^{2\pi}_{4\pi/3}\cos(\xi/2)
\Phi(\xi-2\pi)\Phi(\xi/2)\,d\xi.
\end{equation}
Observe that, for any $\xi\in[4\pi/3,2\pi]$,
by \eqref{ph3}, we know that $\Phi(\xi-2\pi)=1/\sqrt{2\pi}$.
Moreover, by $0\le \frac{\xi}2\frac{3}{2\pi}-1\le\frac12$
and the construction of $g$, we conclude
that $0=g(0)\le g(\frac{\xi}2\frac{3}{2\pi}-1)\le g(\frac12)$.
Hence, from the construction of $\Phi$, we deduce that
$$
1\ge\Phi(\xi/2)\ge\frac1{\sqrt{2\pi}}\cos\lf(\frac{\pi}2
g\lf(\frac12\r)\r),
$$
which, together with \eqref{x.r}, further implies that
\begin{align*}
\lf|{\rm J}_3\r|
&\ge-2\int_{4\pi/3}^{2\pi}\cos(\xi/2)\frac1{\sqrt{2\pi}}
\cos\lf(\frac{\pi}2g\lf(\frac12\r)\r)\frac1{\sqrt{2\pi}}\,d\xi\\
&=-\frac1{\pi}\cos\lf(\frac{\pi}2g\lf(\frac12\r)\r)
\int_{4\pi/3}^{2\pi}\cos(\xi/2)\,d\xi
=\frac{\sqrt{3}}{\pi}\cos\lf(\frac{\pi}2g\lf(\frac12\r)\r).
\end{align*}

Combining the estimates of ${\rm J}_1$, ${\rm J}_2$ and ${\rm J}_3$
and the construction of $g$, we obtain
\begin{align*}
|\psi(0)|&\ge\lf|{\rm J}_3+{\rm J}_2\r|-\lf|{\rm J}_1\r|
=\lf|{\rm J}_3\r|+\lf|{\rm J}_2\r|-\lf|{\rm J}_1\r|\\
&\ge\frac{\sqrt{3}}{\pi}\cos\lf(\frac{\pi}2g\lf(\frac12\r)\r)
+\frac{2-\sqrt{3}}{\pi}\cos\lf(\frac{\pi}2g\lf(\frac12\r)\r)
-\frac{2-\sqrt{3}}{\pi}\cos\lf(\frac{\pi}2g\lf(\frac12\r)\r)\\
&=\frac{\sqrt{3}}{\pi}\cos\lf(\frac{\pi}2g\lf(\frac12\r)\r)>0,
\end{align*}
which completes the proof of Example \ref{ec.a}.
\end{example}

Now we are ready to prove Theorem \ref{ta.x}.

\begin{proof}[Proof of Theorem \ref{ta.x}]
(I) Suppose that $q\in[2,\fz)$, $\phi$ and $\Phi$ are
defined as in the construction of the $1$-dimensional Meyer wavelets.
Moreover, we assume that the $1$-dimensional Meyer wavelet
$\psi$ satisfies $\psi(0)\neq0$.

For any $x:=(x_1,\ldots,x_D)\in\rr^D$, let
$\psi^{\vec{0}}(x):=\phi(x_1)\cdots\phi(x_D)$.
From \cite[(5.1)]{lly}, we deduce that, for all $x\in\rr^D$,
\begin{equation}\label{c.a}
\sum_{k\in\zz^D}\lf|\psi^{\vec{0}}(x-k)\r|
=\prod_{\ell=1}^D\sum_{k_\ell\in\zz}\lf|\phi\lf(x_\ell-k_\ell\r)\r|\ls1.
\end{equation}

For any $j\in\zz$, $k\in\zz^D$ and $x\in\rr^D$, we write
$\psi^{\vec{0}}_{j,\,k}(x):=2^{Dj/2}\psi^{\vec{0}}\lf(2^jx-k\r)$.
Let $f\in\lon$. For any $j\in\zz$, define
$P_j(f):=\sum_{k\in\zz^D}\langle f,\psi^{\vec{0}}_{j,\,k}\rangle
\psi^{\vec{0}}_{j,\,k}$. Then, for any $j\in\zz$, by \eqref{c.a},
\begin{align}\label{c.b}
\lf\|P_j(f)\r\|_{\lon}&\le\int_{\rr^D}\int_{\rr^D}
|f(y)|\sum_{k\in\zz^D}\lf|\psi^{\vec{0}}(2^jy-k)\r|
\lf|2^{Dj}\psi^{\vec{0}}(2^jx-k)\r|\,dx\,dy\\
&\noz\ls\int_{\rr^D}\int_{\rr^D}
|f(y)|\sum_{k\in\zz^D}\lf|\psi^{\vec{0}}(x-k)\r|\,dy\ls\|f\|_{\lon}.
\end{align}

The inclusion relation in (i) of Theorem \ref{ta.x} is easy to see.
In order to prove (i) of Theorem \ref{ta.x}, it suffices to show that
$\tls\subsetneqq  L^{1}( \rr^D )\bigcup \tls$.
We first observe that $\psi^{\vec{0}}\in\lon$. Indeed,
$$
\lf\|\psi^{\vec{0}}\r\|_{\lon}=\prod_{\ell=1}^D\lf\|\phi\r\|_{\lon}
=\lf\|\phi\r\|_{\lon}^D<\fz.
$$

To show $\psi^{\vec{0}}\not\in\tls$, let
$a_{j,\,k}(\psi^{\vec{0}}):=(\psi^{\vec{0}},\psi_{j,\,k})$
for any $j\in\zz$ and $k\in\zz^D$, where
$(\cdot,\cdot)$ represents the $\ltw$ inner product.
Let $\xi:=(\xi_1,\ldots,\xi_D),\,\eta:=(\eta_1,\ldots,\eta_D)\in\rr^D$.
Then, by the multiplication formula (see \cite[p.\,8, Theorem 1.15]{sw}),
\eqref{ph4}, \eqref{ph2}, \eqref{ps1}
and the assumption $\psi(0)\neq0$, we obtain
\begin{align}\label{3.4x}
\lf|a_{j,\,0}\lf(\psi^{\vec{0}}\r)\r|&=\prod_{\ell=1}^D\lf|\int_{\rr}
\widehat{\phi}\lf(-\xi_\ell\r)
2^{-j/2}\widehat{\psi}\lf(2^{-j}\xi_\ell\r)\,d\xi_\ell\r|\\
&=\prod_{\ell=1}^D\lf|\int_{-4\pi/3}^{4\pi/3}
\Phi\lf(\xi_\ell\r)
2^{-j/2}\widehat{\psi}\lf(2^{-j}\xi_\ell\r)\,d\xi_\ell\r|\noz\\
&=\prod_{\ell=1}^D\lf|\int_{-2^{2-j}\pi/3}^{2^{2-j}\pi/3}
\Phi\lf(2^j\eta_\ell\r)2^{j/2}\widehat{\psi}\lf(\eta_\ell\r)\,d\eta_\ell\r|\noz\\
&=\prod_{\ell=1}^D\lf|\int_{-8\pi/3}^{8\pi/3}
\Phi\lf(2^j\eta_\ell\r)2^{j/2}\widehat{\psi}\lf(\eta_\ell\r)\,d\eta_\ell\r|\noz\\
&\sim2^{Dj/2}\prod_{\ell=1}^D\lf|\int_{-8\pi/3}^{8\pi/3}
\widehat{\psi}\lf(\eta_\ell\r)
\,d\eta_\ell\r|\sim2^{Dj/2}|\psi(0)|^D\gtrsim2^{Dj/2}\noz,
\end{align}
provided that $j<-M$ for some positive integer $M$ large enough.
Therefore, we have
\begin{align*}
&\int_{\rr^D}\int_{\rr^D}
\lf\{\sum_{j\in\zz,\,k\in\zz^D}\lf[2^{Dj/2}\lf|a_{j,\,k}\lf(\psi^{\vec{0}}\r)\r|
\chi\lf(2^jx-k\r)\r]^q\r\}^{1/q}\,dx\\
&\hs\ge\int_{\rr^D}
\lf\{\sum_{j=-\fz}^{-M-1}2^{Djq/2}\lf|a_{j,\,0}\lf(\psi^{\vec{0}}\r)\r|^q
\chi\lf(2^jx-k\r)\r\}^{1/q}\,dx\\
&\hs\gtrsim\int_{\rr^D}
\lf\{\sum_ {j=-\fz}^{-M-1}2^{Djq}\chi\lf(2^jx-k\r)\r\}^{1/q}\,dx\\
&\hs\gtrsim\sum_{m=M}^\fz\int_{B(0,\,2^{m+1})\bh B(0,\,2^m)}
\lf\{\sum_{j=-\fz}^{-m-1}2^{Djq}\r\}^{1/q}\,dx=\fz,
\end{align*}
which, combined with Theorem \ref{ta.d}, implies that
$\psi^{\vec{0}}\not\in\tls$.
This finishes the proof of (i) of Theorem \ref{ta.x}.

(II) Then we use Daubechies wavelets to proof  (ii) of Theorem \ref{ta.x}.
We know that there exist some integer $M$
and a Daubechies scale function $\Phi^{0}(x)\in C^{D+2}_{0} ([-2^{M}, 2^{M}]^{D})$
satisfying
\begin{equation}\label{5.1}
C_{D}=\int \frac{-y_{1}}{|y|^{n+1}}
\Phi^{0}(y-2^{M+1}e) dy<0,\end{equation}
where $e= (1,1,\cdots, 1)$.
Let $\Phi(x)= \Phi^{0}(x-2^{M+1}e)$ and let
$f$ be defined as
\begin{equation}\label{5.2}
f(x)=\sum\limits_{j\in 2\mathbb{N}} \Phi(2^{j}x).
\end{equation}

For $j,j'\in 2\mathbb{N}, j\neq j'$,  the
supports of $\Phi(2^{j}x)$ and  $\Phi(2^{j'}x)$ are
disjoint. Hence the above $f(x)$ in (\ref{5.2}) belongs to $L^{\infty}(\mathbb{R}^{D})$.
The same reasoning gives, for any $j'\in \mathbb{N}$,
$$\sum\limits_{j\in \mathbb{N}, 2j> j'} \Phi(2^{2j}x)\in
L^{\infty}(\mathbb{R}^{D}).$$ Now we compute the wavelet coefficients of $f(x)$ in (\ref{5.2}).
For $(\lambda',j',k')\in \Lambda_{D}$, let $f^{\lambda'}_{j',k'}=
\langle f,\ \Phi^{\lambda'}_{j',k'}\rangle$. We divide
two cases: $j'<0$ and $j'\geq 0$.

For $j'<0$, since the support of $f$ is contained in $[-3 \cdot
2^{M}, 3 \cdot 2^{M}]^{D}$, we know that if $|k'|> 2^{2M+5}$, then
$f^{\lambda'}_{j',k'}=0$. If $|k'|\leq 2^{2M+5}$, we have
$$|f^{\lambda'}_{j',k'}| \leq C 2^{Dj'} \int |f(x)| dx
\leq C 2^{Dj'}.$$ For $j'\geq 0$, by orthogonality of the
wavelets, we have
$$f^{\lambda'}_{j',k'}= \Big\langle f,\ \Phi^{\lambda'}_{j',k'}\Big\rangle
= \Big\langle \sum\limits_{j\in \mathbb{N}, 2j> j'} \Phi(2^{2j}\cdot),\
\Phi^{\lambda'}_{j',k'}\Big\rangle.$$
By the same reasoning,
for the case $j'\geq0$, we know that if $|k'|> 2^{2M+5}$, then
$f^{\lambda'}_{j',k'}=0$. Since $\sum\limits_{j\in \mathbb{N}, 2j>
j'} \Phi(2^{2j}x)\in L^{\infty}$, if $|k'|\leq 2^{2M+5}$, we
have $$|f^{\lambda'}_{j',k'}| \leq C\int
|\Phi^{\lambda'}_{j',k'}(x)| dx \leq C .$$ By the
above estimation of wavelet coefficients of $f(x)$ and by the wavelet
characterization of $\dot{F}^{0}_{\infty,q'}(\mathbb{R}^{D})$ in (ii) of Theorem \ref{ta.d}, we conclude that $f\in \dot{F}^{0}_{\infty,q'}(\mathbb{R}^{D})$.
That is to say,
\begin{equation}\label{eq:1111}f\in  L^{\infty}(\mathbb{R}^{D})\bigcap \dot{F}^{0}_{\infty,q'}(\mathbb{R}^{D}).\end{equation}

Since $\Phi^{0}\in C^{D+2}_{0}([-2^{M}, 2^{M}]^{D})$, we know that
$$\Phi(x)= \Phi^{0}(x-2^{M+1}e)\in C^{D+2}_{0}([2^{M}, 3 \cdot 2^{M}]^{D}).$$
Further,  if $|x|\leq 2^{M-1}$ and $y\in [2^{M},
3 \cdot 2^{M}]^{D} $, then $|x-y|> 2^{M-1}$. Hence $R_{1}\Phi(x)$ is
smooth in the ball $\{x:\ |x|\leq 2^{M-1}\}$.

Applying (\ref{5.1}), there exists a
positive $\delta>0$ such that for $|x|<\delta$, there holds
$R_{1}\Phi(x)<\frac{C_{D}}{2}<0.$ That is to say, if
$2^{2j}|x|<\delta$, then $R_{1}\Phi(2^{2j}x)<\frac{C_{D}}{2}<0$.
Hence \begin{equation}\label{eq:1112} R_{1}f(x)\notin L^{\infty}(\mathbb{R}^{D}).\end{equation}

The equations (\ref{eq:1111}), (\ref{eq:1112}) and the continuity of Riesz operators on $\dot{F}^{0}_{\infty,q'}(\mathbb{R}^{D})$ implies the conclusion (ii).

\end{proof}

\section{The proof of Theorems \ref{th:111}, \ref{ta.i} and \ref{ta.cor}}
We prove first Theorem \ref{th:111}
\begin{proof}
If $l \in (L^{1}(\mathbb{R}^{D})\bigcup \dot{F}^{0}_{1,q}(\mathbb{R}^{D}))'$, then
$$\sup\limits_{f\in \schi, \|f\|_{L^{1} \bigcup \dot{F}^{0}_{1,q} }\leq 1} |\langle l, f\rangle|<\infty.$$
That is to say,
\begin{equation}\label{e1e}\sup\limits_{f\in \schi, \|f\|_{L^{1} }\leq 1} |\langle l, f\rangle|<\infty \mbox { and } \end{equation}
\begin{equation}\label{e2e}\sup\limits_{f\in \schi, \|f\|_{ \dot{F}^{0}_{1,q} }\leq 1} |\langle l, f\rangle|<\infty.\end{equation}
The condition (\ref{e1e}) means $l\in L^{\infty}(\mathbb{R}^{D})$,
the condition (\ref{e2e}) means $l\in \dot{F}^{0}_{\infty,q'}(\mathbb{R}^{D})$.
Hence we have the following inclusion relation:
\begin{equation}\label{eq:inc.1}
\big(L^{1}(\mathbb{R}^{D})\bigcup \dot{F}^{0}_{1,q}(\mathbb{R}^{D})\big)' \subset L^{\infty}(\mathbb{R}^{D}) \bigcap \dot{F}^{0}_{\infty,q'}(\mathbb{R}^{D}).
\end{equation}

Further, it is known that
$$L^{1}(\mathbb{R}^{D})\bigcup \dot{F}^{0}_{1,q}(\mathbb{R}^{D}) \subset {\rm WE}^{1,q} (\mathbb{R}^{D}) \subset L^{1}(\mathbb{R}^{D})+ \dot{F}^{0}_{1,q}(\mathbb{R}^{D}).$$
Hence we have
\begin{equation}\label{eq:inc.2}
\big(L^{1}(\mathbb{R}^{D})+ \dot{F}^{0}_{1,q}(\mathbb{R}^{D})\big)'\subset \big( {\rm WE}^{1,q} (\mathbb{R}^{D})\big)' \subset \big(L^{1}(\mathbb{R}^{D})\bigcup \dot{F}^{0}_{1,q}(\mathbb{R}^{D})\big)'.
\end{equation}

Moreover, we know that
\begin{equation}\label{eq:inc.3}
\big(L^{1}(\mathbb{R}^{D})+ \dot{F}^{0}_{1,q}(\mathbb{R}^{D})\big)'= L^{\infty}(\mathbb{R}^{D}) \bigcap \dot{F}^{0}_{\infty,q'}(\mathbb{R}^{D}).
\end{equation}

The equations (\ref{eq:inc.1}) , (\ref{eq:inc.2}) and (\ref{eq:inc.3}) implies the Theorem \ref{th:111}.

\end{proof}

Then we prove Theorem \ref{ta.i}.
\begin{proof}
By the continuity of the Riesz operators on the $\dot{F}^{0}_{\infty,q}(\mathbb{R}^{D})$, we know that if  $f_{l}\in \dot{F}^{0}_{\infty,q}(\mathbb{R}^{D})\bigcap
L^{\infty}(\mathbb{R}^{D})$, then
$$\sum\limits_{0\leq l\leq D}R_{l}f_{l}(x) \in
\dot{F}^{0}_{\infty,q}(\mathbb{R}^{D}).$$

Now we prove the converse result. Let
$$\begin{array}{c} B=\Big\{(g_{0},g_{1},\cdots,g_{D}): g_{l}\in {\rm WE}^{1,q'}(\mathbb{R}^{D}),
l=0,\cdots,D\Big\},\\
\tilde{B}=\Big\{(g_{0},g_{1},\cdots,g_{D}): g_{l}\in L^{1}(\mathbb{R}^{D}) + \dot{F}^{0}_{1,q'}(\mathbb{R}^{D}),
l=0,\cdots,D\Big\},\end{array}$$
where $B\subset \tilde{B}$. The norm of $B$ and $\tilde{B}$ is defined as follows respectively $$ \begin{array}{c}\|(g_{0},g_{1},\cdots,g_{D})\|_{B}=
\sum\limits^{D}_{l=0}\|g_{l}\|_{{\rm WE}^{1,q'}},\\
\|(g_{0},g_{1},\cdots,g_{D})\|_{\tilde{B}}=
\sum\limits^{D}_{l=0}\|g_{l}\|_{L^{1} + \dot{F}^{0}_{1,q'}}.\end{array}$$
 We define
$$\begin{array}{c}S=\Big\{(g_{0},g_{1},\cdots,g_{D})\in B: g_{l}=R_{l}g_{0}, l=0, 1,\cdots, D\Big\},\\
\tilde{S}=\Big\{(g_{0},g_{1},\cdots,g_{D})\in \tilde{B}: g_{l}=R_{l}g_{0}, l=0, 1,\cdots, D\Big\},\end{array}$$
where $S\subset \tilde{S}.$

By Theorem \ref{ta.h}, $g_{0}\rightarrow
(g_{0},R_{1}g_{0},\cdots,R_{D}g_{0})$ define a norm preserving map
from $\dot{F}^{0}_{1,q'}(\mathbb{R}^{D})$ to $S$. Hence the set of continuous linear functionals $f$ on
$\dot{F}^{0}_{1,q'}(\mathbb{R}^{D})$ is equivalent to the set of bounded linear map on the set $S$.
According to Theorem \ref{th:111}, which is also the
set of continuous linear functionals on Banach space $\tilde{S}$.

According to Theorem \ref{th:111}, the continuous linear functionals on $B$ belong to
$${\rm WE}^{\infty,q}(\mathbb{R}^{D}) + \cdots + {\rm WE}^{\infty,q}(\mathbb{R}^{D}).$$
$\forall f\in \dot{F}^{0}_{\infty,q}(\mathbb{R}^{D})$, $f$ defines a continuous linear functional
$l$ on  $\dot{F}^{0}_{1,q'}(\mathbb{R}^{D})$ and also on $\tilde{S}$. Hence there exist
$\tilde{f}_{l}\in {\rm WE}^{\infty,q}(\mathbb{R}^{D}), l=0,1,\cdots,D,$ such that for any $g_{0}\in \dot{F}^{0}_{1,q'}(\mathbb{R}^{D})$,
\begin{eqnarray*}
 &  &\int_{\mathbb{R}^{D}}f(x) g_{0}(x) dx \\
 &=&\int_{\mathbb{R}^{D}}\tilde{f}_{0}(x) g_{0}(x) dx
+\sum\limits^{D}_{l=1} \int_{\mathbb{R}^{D}}
\tilde{f}_{l}(x) R_{l}g_{0}(x) dx\\
&=&\int_{\mathbb{R}^{D}}\tilde{f}_{0}(x) g_{0}(x) dx
-\sum\limits^{D}_{l=1} \int_{\mathbb{R}^{D}} R_{l}(\tilde{f}_{l})(x)
g_{0}(x) dx.
\end{eqnarray*}

Hence $f(x)= \tilde{f}_{0}(x) -\sum\limits^{D}_{l=1}
R_{l}(\tilde{f}_{l})(x)$.

\end{proof}

Finally, we give the proof of Theorem \ref{ta.cor}.
\begin{proof}
By the continuity of Riesz operators on   $  \dot{F}^{0}_{1,q}(\mathbb{R}^{D}$, there exists a positive constant $C$
such that, for all $f\in\tls$,
$$\sum_{\ell=0}^D\|R_{\ell}(f)\|_{L^{1}(\mathbb{R}^{D}) + \dot{F}^{0}_{1,q}(\mathbb{R}^{D})}
\le C\|f\|_{\tls}.
$$

To prove
$$
\frac1C\|f\|_{\tls}\le\sum_{\ell=0}^D\|R_{\ell}(f)\|_{L^{1}(\mathbb{R}^{D}) + \dot{F}^{0}_{1,q}(\mathbb{R}^{D})},
$$
it is sufficient to prove
$$|\langle f,g\rangle|\leq C\big\{\sum\limits^{D}_{l=0} \|R_l f\|_{{\rm WE}^{1,q}(\mathbb{R}^{D})}\big\}\|g\|_{\dot{F}^{0}_{\infty,q'}(\mathbb{R}^{D})},
\forall g\in \schi \bigcap \dot{F}^{0}_{\infty,q'}(\mathbb{R}^{D}).$$
But $\forall g\in \schi \bigcap\dot{F}^{0}_{\infty,q'}(\mathbb{R}^{D})$, by Theorem \ref{ta.i},  there exists $g_{l}$ such that $\|g_l\|_{{\rm WE}^{\infty,q'}(\mathbb{R}^{D})}\leq C\|g\|_{\dot{F}^{0}_{\infty,q'}(\mathbb{R}^{D})}$ and $g=\sum\limits^{D}_{l=0} R_l g_l.$
Hence, we have
$$\begin{array}{rcl}
|\langle f, g\rangle | &=  |\langle f,  \sum\limits^{D}_{l=0} R_l g_l\rangle |\,\,\, \leq \sum\limits^{D}_{l=0}  |\langle f,  R_l g_l\rangle |
&= \sum\limits^{D}_{l=0}  |\langle R_lf,   g_l\rangle | \\ &\leq C \sum\limits^{D}_{l=0} \|R_l f\|_{{\rm WE}^{1,q}(\mathbb{R}^{D})} \|g_{l}\|_{{\rm WE}^{\infty,q'}(\mathbb{R}^{D})}
& \leq C \sum\limits^{D}_{l=0} \|R_l f\|_{{\rm WE}^{1,q}(\mathbb{R}^{D})} \|g\|_{\dot{F}^{0}_{\infty,q'}(\mathbb{R}^{D})}.
\end{array}$$
\end{proof}

{\bf Acknowledgement:}
The authors would like to thank Dachun Yang and Xing Fu, who contributed beneficial discussions and useful suggestions to this study.

\bigskip

\noindent Qixiang Yang

\medskip

\noindent School of Mathematics and Statistics,
Wuhan University, Wuhan, 430072, China.

\smallskip

\noindent{\it E-mail:} \texttt{qxyang@whu.edu.cn}

\bigskip

\noindent  Tao Qian (Corresponding author)

\medskip

\noindent  Department of Mathematics, University of Macau, Macao, China

\smallskip

\noindent {\it E-mails}: \texttt{ fsttq@umac.mo} 

\end{document}